\documentclass[12pt]{article}
\usepackage{amssymb,amsmath,amsthm}
\usepackage{graphicx,epstopdf}
\usepackage{cite,color}
\usepackage{float}
\usepackage{tabularx}
\usepackage{graphicx}
\newtheorem{theorem}{Theorem}[section]

\newtheorem{lemma}{Lemma}[section]
\newtheorem{definition}{Definition}[section]
\newtheorem{remark}{Remark}[section]

\newcommand{\beq}{\begin{equation}}
	\newcommand{\eeq}{\end{equation}}
\newcommand{\beqn}{\begin{eqnarray}}
	\newcommand{\eeqn}{\end{eqnarray}}

\def\to{{\rightarrow}}

\allowdisplaybreaks \baselineskip 20pt

\title{Besicovitch Pseudo Almost Periodic Stochastic Processes and Their Applications to Delayed Stochastic Neural Networks\thanks{This work is supported by the National Natural Science Foundation of China
under Grant No. 12261098.}}
\author{ Yongkun Li$^1$\thanks{The corresponding author. Email: yklie@ynu.edu.cn.}, Yuwei Cao$^1$ and Bing Li$^2$\\
$^1$Department of Mathematics,\\
Yunnan University, Kunming, Yunnan 650091, China\\
	$^2$School of Mathematics and Computer Science,\\
Yunnan Minzu University, Kunming, Yunnan 650500,  China\\
}
\date{}
\begin{document}
\maketitle{}	
\begin{abstract}
The Besicovitch almost periodic concept is a generalization of the Bohr, Stepanov, and Weyl almost periodic concepts and, in a certain sense, has broader coverage and weaker conditions. However, the concept of pseudo almost periodic stochastic processes in the Besicovitch sense has not yet been systematically proposed, which limits the characterization of the complex oscillatory behavior of a broader class of stochastic non-autonomous systems. This paper proposes the concept of Besicovitch pseudo almost periodic stochastic processes. As an application of this concept, we study the existence of Besicovitch pseudo almost periodic solutions in distribution for a class of stochastic delayed neural networks  in which all parameters except the time delays are pseudo almost periodic and establish stability criteria for such solutions in the $p$-th mean sense. The work of this paper fills the gap caused by the lack of the concept of Besicovitch pseudo almost periodic stochastic processes and provides a new conceptual tool and analytical approach for characterizing the complex oscillatory behavior of stochastic differential equations and stochastic neural network systems.
\end{abstract}
	
	\allowdisplaybreaks
	\section{Introduction}
	\setcounter{equation}{0}
	\indent

Since Bohr introduced the concept of almost periodic functions \cite{bo}, the theory of almost periodicity has been continuously generalized and developed, with concepts such as Stepanov almost periodicity \cite{st}, Weyl almost periodicity \cite{we}, and Besicovitch almost periodicity \cite{be}, as well as pseudo almost periodicity in various senses, appearing successively \cite{y1,wave}. The concept of almost periodicity and its variants not only expand the scope of describing periodic phenomena but also provide an effective tool for characterizing the long-term behavior of non-autonomous systems with complex oscillatory properties \cite{wave,wave1}. Random phenomena are ubiquitous in natural and social systems. Consequently, concepts of almost periodic stochastic processes in various senses have also been proposed successively, and studying the existence and stability of almost periodic solutions of stochastic differential equations (SDEs) in various senses has become one of the important problems in the theory and applications of SDEs and dynamical systems.

In the stochastic framework, the Besicovitch almost periodic concept is a generalization of the Bohr, Stepanov, and Weyl almost periodic concepts and, in a certain sense, has broader coverage and weaker conditions \cite{y2}. Existing work has proposed concepts of Besicovitch almost periodic stochastic processes in the $p$-th mean and distribution senses, respectively, and, taking stochastic neural networks (NNs) as an example, discussed the existence and stability of almost periodic solutions in distribution \cite{y2}. However, as far as the authors know, the concept of pseudo almost periodic stochastic processes in the Besicovitch sense has not yet been systematically proposed. This conceptual gap limits the characterization of the complex oscillatory behavior of a broader class of stochastic non-autonomous systems. Moreover, for general SDEs, counterexamples have shown that nontrivial solutions of some SDEs with almost periodic coefficients may not be almost periodic in the mean square sense but may be almost periodic in distribution; this means that the requirement of almost periodic solutions in the $p$-th mean sense is too strong, and establishing the corresponding theory of almost periodic solutions in distribution has clear necessity and rationality \cite{y3}.

On the one hand, the connection weights and external inputs of NNs usually vary with time, and the differential equations describing NNs are usually non-autonomous. For non-autonomous systems, constant equilibrium solutions usually do not exist; even if all coefficients are periodic functions, such NNs usually have no periodic solutions if their periods are incommensurable. Therefore, for non-autonomous NNs, the existence of almost periodic solutions is as important as the existence of equilibrium solutions for autonomous NNs. Stability is a basic starting point for the design, implementation, and application of NNs, and stability in $p$-th mean is one of the most basic and tractable stability concepts in the analysis of stochastic NNs \cite{y5,y6}. Studying the existence and the stability in $p$-th mean of almost periodic solutions of stochastic NNs has important theoretical and practical significance. In recent years, Besicovitch almost periodic solutions of deterministic NNs have begun to attract attention, for example, the existence and finite-time stability for fractional-order quaternion-valued NNs with time-varying delays \cite{y4}; however, Besicovitch pseudo almost periodic solutions under stochastic perturbations and their stability in $p$-th mean have not yet been systematically treated.

On the other hand, it should be pointed out that studying Besicovitch pseudo almost periodic solutions of stochastic NNs is only one application of the new concept proposed in this paper, rather than the sole purpose of the theoretical construction in this paper. The core intention of this paper is to establish a general concept of pseudo almost periodic stochastic processes in the Besicovitch sense, provide a more general descriptive tool for the complex oscillatory behavior of stochastic non-autonomous systems than the existing Besicovitch almost periodic framework in distribution, and then take a class of stochastic NNs with delays as an application vehicle for this concept to test its operability.

Inspired by the above discussion, this paper aims to propose the concept of Besicovitch pseudo almost periodic stochastic processes and give the corresponding definition of almost periodic solutions in distribution; as an application, we study the existence of Besicovitch pseudo almost periodic solutions in distribution for a class of stochastic delayed NNs in which all parameters except the time delays are pseudo almost periodic and establish stability criteria for such solutions in the $p$-th mean sense. The work of this paper, on the one hand, fills the gap caused by the lack of the concept of Besicovitch pseudo almost periodic stochastic processes; on the other hand, it provides a new conceptual tool and analytical approach for characterizing the complex oscillatory behavior of SDEs and stochastic neural network systems.

The paper is organized as follows. Section 2 introduces some basic concepts and lemmas and describes the stochastic NNs considered in this paper. Section 3 presents the concept of Besicovitch pseudo almost periodic stochastic processes. Section 4 establishes the existence and uniqueness of Besicovitch pseudo almost periodic solutions in distribution for the stochastic NNs considered in this paper, and proves their global exponential stability in $p$-th mean. Section 5 provides a numerical example and simulations to illustrate the effectiveness of the results obtained in this paper. Finally, Section 6 gives the conclusions.
	
\section{Preliminaries and model description  }
\setcounter{equation}{0}
\indent
	
Throughout, let $(\mathbb{X},\|\cdot\|)$ be a separable Banach space. For $1\le p<\infty$, the symbol $\mathcal{L}_{\mathrm{loc}}^p(\mathbb{R},\mathbb{X})$ stands for the set formed by all maps $f:\mathbb{R}\to \mathbb{X}$ whose restriction to any compact interval is $p$-integrable.

We write $\mathcal{L}^\infty(\mathbb{R},\mathbb{X})$ for the collection of all measurable functions $g:\mathbb{R}\to \mathbb{X}$ such that the essential supremum
$\|g\|_\infty:=\mathrm{ess}\sup_{t\in \mathbb{R}}\|g(t)\|$
is finite. Endowed with the norm $\|\cdot\|_\infty$, the space $(\mathcal{L}^\infty(\mathbb{R},\mathbb{X}),\|\cdot\|_\infty)$ becomes a Banach space.

Given $g \in \mathcal{L}_{\mathrm{loc}}^p(\mathbb{R},\mathbb{X})$, define its Besicovitch seminorm by
\[
\|g\|_{B^p}
=\left(\overline{\lim_{l \rightarrow \infty}}
\frac{1}{2l}\int_{-l}^{l}\|g(t)\|^p\,\mathrm{d}t\right)^{1/p}.
\]

\begin{definition}
We say that a function $f \in \mathcal{L}_{\mathrm{loc}}^p(\mathbb{R},\mathbb{X})$ is $B^p$-bounded provided that $\|f\|_{B^p}<\infty$. The family consisting of all such $B^p$-bounded functions is denoted by $B_B^p(\mathbb{R},\mathbb{X})$.
\end{definition}

\begin{definition}
A function $f \in \mathcal{L}_{\mathrm{loc}}^p(\mathbb{R},\mathbb{X})$ is called $B^p$-continuous whenever
$\lim_{\xi \rightarrow 0}\|f(\cdot+\xi)-f(\cdot)\|_{B^p}=0.$
The collection of all $B^p$-continuous functions is written as $B_C^p(\mathbb{R},\mathbb{X})$.
\end{definition}

The common subclass consisting of functions that are simultaneously $B^p$-bounded and $B^p$-continuous will be denoted by $B_{BC}^p(\mathbb{R},\mathbb{X})$.

Let $(\Omega,\mathcal{F},\{\mathcal{F}_t\}_{t\in\mathbb{R}},P)$ be a complete filtered probability space, whose filtration $\{\mathcal{F}_t\}_{t\in\mathbb{R}}$ is natural and fulfills the usual conditions. By $\mathcal{L}^p(\Omega,\mathbb{X})$ we mean the space of $\mathbb{X}$-valued random variables possessing finite $p$-th moment, that is,
$E\|X\|^p=\int_{\Omega}\|X\|^p\,\mathrm{d}P<\infty.$
When supplied with the norm $\|X\|_p:=(E\|X\|^p)^{1/p}$, this space is a Banach space.

Let $\mathcal{L}_{\mathrm{loc}}^p(\mathbb{R},\mathcal{L}^p(\Omega,\mathbb{X}))$ consist of all stochastic processes $X:\mathbb{R}\to \mathcal{L}^p(\Omega,\mathbb{X})$ that are locally $p$-integrable. For $X \in \mathcal{L}_{\mathrm{loc}}^p(\mathbb{R},\mathcal{L}^p(\Omega,\mathbb{X}))$, its Besicovitch seminorm is defined by
\[
\|X\|_{SB^p}
=\left(\overline{\lim_{l \rightarrow \infty}}
\frac{1}{2l}\int_{-l}^{l}E\|X(t)\|^p\,\mathrm{d}t\right)^{1/p}.
\]

\begin{definition}
A stochastic process $X \in \mathcal{L}_{\mathrm{loc}}^p(\mathbb{R},\mathcal{L}^p(\Omega,\mathbb{X}))$ is termed $SB^p$-bounded if $\|X\|_{SB^p}<\infty$. The family of all such $SB^p$-bounded stochastic processes is denoted by $SB_B^p(\mathbb{R},\mathcal{L}^p(\Omega,\mathbb{X}))$.
\end{definition}

\begin{definition}
A stochastic process $X \in \mathcal{L}_{\mathrm{loc}}^p(\mathbb{R},\mathcal{L}^p(\Omega,\mathbb{X}))$ is said to be $SB^p$-continuous whenever
\[
\lim_{\xi \to 0}\|X(\cdot+\xi)-X(\cdot)\|_{SB^p}=0.
\]
We use $SB_C^p(\mathbb{R},\mathcal{L}^p(\Omega,\mathbb{X}))$ to denote the set of all $SB^p$-continuous stochastic processes.
\end{definition}

The class of processes that possess both properties, namely $SB^p$-boundedness and $SB^p$-continuity, is written as $SB_{BC}^p(\mathbb{R},\mathcal{L}^p(\Omega,\mathbb{X}))$. In analogy with the argument for Theorem 2.1 in \cite{2026}, one verifies without difficulty the following assertion.

\begin{lemma}
With respect to the seminorm $\|\cdot\|_{SB^p}$, the space $SB_{BC}^p(\mathbb{R},\mathcal{L}^p(\Omega,\mathbb{X}))$ is complete.
\end{lemma}

    \begin{definition}\label{bp}\cite{2026}
	A function $f \in B_{BC}^p(\mathbb{R},\mathbb{X})$ is said to be $p$-th Besicovitch almost periodic if for every $\varepsilon >0$, there exists a positive number $\ell=\ell(\varepsilon)$ such that every interval of length $\ell(\varepsilon)$ contains a number $\tau$ satisfying
	$\|f(\cdot+\tau)-f(\cdot)\|_{B^p}<\varepsilon.$
	The collection of all such processes will be denoted by $B_{AP}^p(\mathbb{R},\mathbb{X})$.
    \end{definition}

\begin{definition}\label{d1}\cite{W}
	A stochastic process $X \in SB_{BC}^p(\mathbb{R},\mathcal{L}^p(\Omega,\mathbb{X}))$ is said to be Besicovitch almost periodic in $p$-th mean if for every $\varepsilon >0$, there exists a positive number $\ell=\ell(\varepsilon)$ such that every interval of length $\ell(\varepsilon)$ contains a number $\tau$ satisfying
	$\|X(\cdot+\tau)-X(\cdot)\|_{SB^p}<\varepsilon.$
	The collection of all such processes will be denoted by $SB_{AP}^p(\mathbb{R},\mathcal{L}^p(\Omega,\mathbb{X}))$.
\end{definition}

Suppose that $(\mathbb{U},d)$ is a separable complete metric space. Denote by $\mathcal{B}(\mathbb{U})$ the Borel $\sigma$-algebra on $\mathbb{U}$, and by $\mathcal{P}(\mathbb{U})$ the collection of all probability measures defined on $\mathcal{B}(\mathbb{U})$. Let $BC(\mathbb{U})$ be the space consisting of all bounded continuous mappings $\phi:\mathbb{U}\to\mathbb{X}$, endowed with the bounded norm
$\|\phi\|_B:=\sup_{x\in\mathbb{U}}\|\phi(x)\|<\infty.$

For any $\phi\in BC(\mathbb{U})$ and any $\varsigma_1,\varsigma_2\in\mathcal{P}(\mathbb{U})$, we introduce the Lipschitz seminorm, the BL norm, and the BL distance as follows:
\[\|\phi \|_L = \mathop {\sup }\limits_{a \ne b} \frac{\|\phi(a) - \phi(b) \|}{d( a,b)},\,\,\| \phi \|_{BL} = \max \{ \| \phi \|_B ,\| \phi \|_L \},\,\,
	d_{BL}(\varsigma_1,\varsigma_2 ) = \mathop {\sup }\limits_{\|\phi\|_{BL} \le 1} \bigg| \int_\mathbb{U} g d( \varsigma_1-\varsigma_2  )  \bigg|.\]

	\begin{definition}\label{bp1}\cite{W}
		A stochastic process $X:\mathbb{R} \to \mathcal{L}^p(\Omega,\mathbb{X})$ is called p-th Besicovitch almost periodic in distribution if for every $\varepsilon >0$, there exists a positive number $\ell=\ell(\varepsilon)$ such that every interval of length $\ell(\varepsilon)$ contains a number $\tau$ satisfying
		$$\bigg(\overline{\lim_{l \rightarrow \infty}} \frac{1}{2l}\int_{-l}^{l} d_{BL}^p(P \circ [X(t+\tau)]^{-1}, P \circ [X(t)]^{-1}) \mathrm{d} t \bigg)^{\frac{1}{p}}<\varepsilon, $$
		where $P \circ [X(t)]^{-1}$ is the law of $X(t)$ under probability $P$. The collection of such processes will be denoted by $DB_{AP}^p(\mathbb{R},\mathcal{L}^p(\Omega,\mathbb{X}))$.
	\end{definition}

	\begin{lemma}\label{le3.2}\cite{W}
		If a stochastic process $X$ is  Besicovitch almost periodic in $p$-th mean, then it is $p$-th Besicovitch almost periodic in distribution.
	\end{lemma}

\begin{lemma}\cite{a3}(BDG inequality)
	 If $f \in \mathcal{L}^2([t_0,T],\mathbb{R}), p\ge 2, B(t)$ is a Brownian motion, then
	$$
	E\bigg(\sup_{t \in \mathbb{R}}\bigg|\int_{t_0}^{t}f(s) \mathrm{d} B(s)\bigg|\bigg)^p\le C_pE\bigg[\int_{t_0}^{T}|f(s)|^2 \mathrm{d} s \bigg]^{\frac{p}{2}},$$
	where $C_p=\Big(\frac{p^{p+1}}{2(p-1)^{p-1}}\Big)^{\frac{p}{2}}$.
\end{lemma}

\begin{lemma}\cite{27}\label{G}
Let $\varphi:\mathbb{R}\to\mathbb{R}$ be a continuous function satisfying
\begin{align*}
0\le \varphi(t)\le a(t)+b\int_{-\infty}^{t}e^{-c(t-s)}\varphi(s)\,\mathrm{d}s,
\end{align*}
where $a$ is locally integrable, $b,c>0$ with $c>b$, and the above integral is assumed to be convergent. Then, for every $d\in(0,c-d]$ such that
$\int_{-\infty}^{0}e^{ds}a(s)\,\mathrm{d}s<\infty,$
one has
\begin{align*}
\varphi(t)\le a(t)+b\int_{-\infty}^{t}e^{-d(t-s)}a(s)\,\mathrm{d}s.
\end{align*}
\end{lemma}

This paper focuses on the stochastic neural network described below:
\begin{align}\label{e1}
	\mathrm{d} x_i(t)
	=& \bigg[-k_i(t)x_i(t)+\sum_{j = 1}^{n} c_{ij}(t)f_j(x_j(t))+\sum_{j = 1}^{n} d_{ij}(t)g_j(x_j(t-\eta_{ij}(t)))+I_i(t)\bigg]\mathrm{d} t\nonumber\\
	& +\sum_{j = 1}^{n} U_{ij}(t) \delta_j(x_j(t-r_{ij}(t)))\mathrm{d} w_j(t),
\end{align}
where $x_i(t)\in\mathbb{R}$ represents the state of the $i$-th unit at time $t$,  $i=1,2,\ldots,n$; $k_i(t)\in\mathbb{R}$ denotes the self-feedback coefficient at time $t$; $c_{ij}(t)$, $d_{ij}(t)$ and $U_{ij}(t)\in\mathbb{R}$ are connection weights; $f_j,g_j:\mathbb{R}\to\mathbb{R}$ denote activation functions; $I_i(t)\in\mathbb{R}$ stands for the external input received by the $i$-th unit; $\eta_{ij},r_{ij}\in\mathbb{R}^+$ are transmission delays; $w(t)=(w_1(t),w_2(t),\ldots,w_n(t))$, $t\in\mathbb{R}$, is a two-sided $n$-dimensional standard Brownian motion defined on a complete filtered probability space satisfying the usual conditions; and $\delta_j:\mathbb{R}\to\mathbb{R}$ is a Borel measurable function.

Let $BC_{\mathcal{F}0}([-\rho,0],\mathbb{R})$ denote the collection of all bounded, $\mathcal{F}_0$-measurable random variables taking values in $C([-\rho,0],\mathbb{R})$.

The initial data for system \eqref{e1} are prescribed by
\[
x_i(s)=\phi_i(s),\quad s\in[-\rho,0],
\]
where $\phi\in BC_{\mathcal{F}0}([-\rho,0],\mathbb{R})$ and
$\rho=\max_{1\le i,j\le n}\{\sup_{t\in\mathbb{R}}\eta_{ij}(t),\ \sup_{t\in\mathbb{R}}r_{ij}(t)\}.$

\section{Besicovitch pseudo almost periodic stochastic processes}
\setcounter{equation}{0}
	\indent

	This section will introduce the concept of Besicovitch pseudo almost periodic stochastic processes in both the $p$-th mean sense and the distributional sense.

	Let $\Sigma $ indicate the Lebesgue $\sigma $-field of $\mathbb{R}$ and $\mathcal{M}_p$ be the set of all positive measures $\mu $ on $\Sigma $ meeting $\mu (\mathbb{R})=+\infty$ and $\mu ([a,b])<+\infty$ for all $a \le b$.
\begin{definition}\cite{measure1}
	Two measures $\mu$, $\nu \in \mathcal{M}_p$ are said to be equivalent, denoted as $\mu \sim \nu$, if and only if there exist constants $\alpha$, $\beta >0$ and a bounded interval $A$ (eventually $\emptyset$) such that $\alpha \nu (C) \le \mu (C) \le \beta \nu (C)$ for all $C \in \Sigma$ with $C \cap A =\emptyset$.
\end{definition}	
For $\mu \in \mathcal{M}_p$ and $\sigma \in \mathbb{R}$, we define the measure $\mu_{\sigma}$ by $\mu_{\sigma}(A)=\mu(\{a+\sigma,a\in A\})$ for $A\in \Sigma$.

Two basic assumptions:
	\begin{itemize}
	\item [$(A_1)$] Let $\mu$, $\nu \in \mathcal{M}_p$ meet $\overline \lim_{h \rightarrow \infty}  \frac{\mu([-h,h])}{\nu([-h,h])} <\infty$.
	\item [$(A_2)$]For $\mu \in \mathcal{M}_p$ and all $\eta \in \mathbb{R}$, there are $\zeta >0$ and a bounded interval $\mathcal{J} \subset \mathbb{R}$, such that
	$$\mu (a+\eta:a\in A) \le \zeta \mu (A),$$
	provided $A \in \Sigma$ meets $A \cap \mathcal{J}= \emptyset $.	
	\end{itemize}
\begin{lemma}\label{a1}\cite{measure1}
	Let $\mu \in \mathcal{M}_p$, then $\mu$ satisfies $(A_2)$ if and only if $\mu \sim \mu_\sigma$ for all $\sigma \in \mathbb{R}$.
\end{lemma}

\begin{lemma}\cite{measure1}\label{2.5}
Condition $(A_2)$ implies that, for all $\eta > 0$,
\[
\overline{\lim_{h\to \infty}} \frac{\mu([-h - \eta, h + \eta])}{\mu([-h, h])} < \infty.
\]
\end{lemma}
Based on Lemma \ref{2.5}, one can easily show that
\begin{lemma}\label{lem27}
Let $\mu \in \mathcal{M}_{p}$ satisfy condition $(A_2)$. Then, for every $\sigma \in \mathbb{R}$ and $\eta > 0$,
\[
\overline{\lim_{h\to \infty}} \frac{\mu_\sigma([-h-\eta, h+\eta])}{\mu([-h, h])} < \infty.
\]
\end{lemma}

For $\mu$, $\nu \in \mathcal{M}_p$, set
\begin{gather*}
\mathcal{O}^p (\mathbb{R},\mathbb{X},\mu,\nu)=\left\{f \in \mathcal{L}^p_{\mathrm{loc}}(\mathbb{R}, \mathbb{X}):\lim_{h \rightarrow \infty}\frac{1}{\nu ([-h,h])}\int_{-h}^{h}\|f(t)\|^p \mathrm{d} \mu (t)=0\right\},\\
\mathcal{O}_M^p (\mathbb{R},\mathcal{L}^p(\Omega,\mathbb{X}),\mu,\nu)=\left\{X \in \mathcal{L}^p_{\mathrm{loc}}( \mathbb{R},\mathcal{L}^p(\Omega,\mathbb{X})):\lim_{h \rightarrow \infty}\frac{1}{\nu ([-h,h])}\int_{-h}^{h}E\|X(t)\|^p \mathrm{d} \mu (t)=0\right\}.
\end{gather*}

\begin{definition}
	Let $\mu,\nu\in \mathcal{M}_p$. A stochastic process $X\in \mathcal{L}^p_{\mathrm{loc}}(\mathbb{R},\mathcal{L}^p(\Omega,\mathbb{X}))$ belongs to $\mathcal{O}_D^p(\mathbb{R},
	\mathcal{L}^p(\Omega,\mathbb{X}),\mu,\nu)$ if and only if it satisfies
	$$\lim\limits_{h\rightarrow \infty}\frac{1}{\nu([-h,h])}\int_{-h}^{h} d_{BL}^p(P\circ [X(t)]^{-1},\delta_0) \mathrm{d} \mu(t)=0,$$
	where $\delta _{0}$ denotes the distribution of random variable $\mathbf{0}$.
\end{definition}

We are now ready to introduce the notion of Besicovitch pseudo almost periodic stochastic processes.
\begin{definition}\label{def333}
	Let $\mu$, $\nu \in \mathcal{M}_p$. A function $f : \mathbb{R}\rightarrow \mathbb{X}$ is called to be $p$-th Besicovitch $(\mu,\nu )$-pseudo almost periodic if it can be expressed as
	$f=f_1+f_2,$
	where $f_1 \in B_{AP}^p(\mathbb{R},\mathbb{X})$ and $f_2 \in \mathcal{O}^p (\mathbb{R},\mathbb{X},\mu,\nu)$. This class is defined by $B_{PAP}^p(\mathbb{R},\mathbb{X},\mu,\nu )$.
\end{definition}

\begin{definition}\label{d2}
	Let $\mu$, $\nu \in \mathcal{M}_p$. A stochastic process $X: \mathbb{R} \rightarrow \mathcal{L}^p(\Omega,\mathbb{X})$ is called to be  Besicovitch $(\mu,\nu )$-pseudo almost periodic in $p$-th mean if it can be expressed as $X=X_1+X_2,$
	where $X_1 \in SB_{AP}^p(\mathbb{R},\mathcal{L}^p(\Omega,\mathbb{X}))$ and $X_2 \in \mathcal{O}_M^p(\mathbb{R},\mathcal{L}^p(\Omega,\mathbb{X}),\mu,\nu)$. The collection of such processes will be defined by $SB_{PAP}^p(\mathbb{R},\mathcal{L}^p(\Omega,\mathbb{X}),\mu,\nu )$.
\end{definition}

\begin{definition}
		Let $\mu$, $\nu \in \mathcal{M}_p$. A stochastic process $X: \mathbb{R} \rightarrow \mathcal{L}^p(\Omega,\mathbb{X})$ is called to be p-th Besicovitch $(\mu,\nu )$-pseudo almost periodic in distribution if it can be decomposed into $X=X_1+X_2,$
	where $X_1 \in DB_{AP}^p(\mathbb{R},\mathcal{L}^p(\Omega,\mathbb{X}))$ and $X_2 \in\mathcal{O}_D^p (\mathbb{R},\mathcal{L}^p(\Omega,\mathbb{X}),\mu,\nu)$. The collection of such processes will be defined by $DB_{PAP}^p(\mathbb{R},\mathcal{L}^p(\Omega,\mathbb{X}),\mu,\nu)$.
\end{definition}

For a stochastic process $X(t)$, noting the fact that
\begin{align*}
d_{BL}^p(P \circ [X(t)]^{-1},\delta_0) \le \bigg(\int_\Omega  \|X(t)\|\mathrm{d} P \bigg)^p\leq E\|X(t)\|^p.
\end{align*}
We have
\begin{theorem}\label{key}
	If $X\in \mathcal{O}_D^p(\mathbb{R},\mathcal{L}^p(\Omega,\mathbb{X}),\mu,\nu)$, then $X\in \mathcal{O}_M^p(\mathbb{R},\mathcal{L}^p(\Omega,\mathbb{X}),\mu,\nu)$.
\end{theorem}

Then, based on Lemma $\ref{le3.2}$ and Theorem \ref{key}, we can easily get the following result.
\begin{theorem}\label{key1}
	If a stochastic process $X$ is Besicovitch $(\mu,\nu )$-pseudo almost periodic in $p$-th mean, then it is p-th Besicovitch $(\mu,\nu )$-pseudo almost periodic in distribution.
\end{theorem}

 \section{Besicovitch pseudo almost periodic solutions in distribution}
\setcounter{equation}{0}
\indent

In this section, we will investigate the existence and global exponential stability of  Besicovitch $(\mu,\nu )$-pseudo almost periodic solution in distribution for \eqref{e1}.

\begin{definition}\label{keyd1}
	An $\mathcal{F}_t$-progressively measurable stochastic process $x(t)=(x_1(t),x_2(t),\ldots,\\
x_n(t))$ is called a mild solution of system \eqref{e1}, if $x_i(t)$ satisfies the following integral equation
	\begin{align}\label{jie}
		x_i(t) =x_i(t_0)e^{-\int_{t_0}^{t} k_i (u)\mathrm{d} u}+\int_{t_0}^t e^{-\int_s^t k_i (u)\mathrm{d} u}\Theta_i(s,x)\mathrm{d} s+\int_{t_0}^t e^{-\int_s^t k_i (u)\mathrm{d} u}\Xi_{ij} (s,x)\mathrm{d} w_j(s),
	\end{align}
where
\begin{align*}
	\Theta_i(s,x)=&\sum_{j = 1}^{n} c_{ij}(s) f_j(x_j(s))+\sum_{j = 1}^{n} d_{ij}(s) g_j(x_j(s-\eta_{ij}(s)))+I_i(s),\\
	\Xi_{ij} (s,x)=&\sum_{j = 1}^{n}U_{ij}(s) \delta_j(x_j(s-r_{ij}(s))).
\end{align*}
\end{definition}

In the sequel,	for $x=(x_1, x_2, \cdots, x_n) \in \mathbb{R}^n$, we denote $|x|_0=\mathop {\max }\limits_{1 \le i \le n}\{|x_i|\}$. For the sake of convenience, we introduce the following notations:
	$$ k_i^- = \inf_{t\in \mathbb{R}} k_i(t),\,\,  k^- = \min_{1\le i\le n}\{ k_i^-\}, \,\,
  c_{ij}^+ = |c_{ij}|_\infty,\,\,
  d_{ij}^+ = |d_{ij}|_\infty,\,\,  U_{ij}^+ = | U_{ij}|_\infty,\,\, I_i^+= |I_i|_\infty,
  $$
  $$
\eta_{ij}^ +  = \mathop {\sup }\limits_{t \in \mathbb{R}}  \eta_{ij}(t),\,\,
r_{ij}^+= \mathop {\sup }\limits_{t \in \mathbb{R}} r_{ij}(t),\,\, \underline{k}_i^1 = \inf_{t\in \mathbb{R}}| k_i^1(t) |, \,\, \bar c_{ij}^1=|c_{ij}^1|_\infty, \,\, \bar d_{ij}^1=|d_{ij}^1|_\infty,
\,\, \bar U_{ij}^1=|U_{ij}^1|_\infty,
$$
and make the following hypothesises:
\begin{itemize}
	\item [$(H_1)$]
	For $i$, $j=1,2, \cdots ,n$, $k_i:\mathbb{R} \rightarrow  \mathbb{R}^+$ with $k_i^-\geq \underline{k}_i^1>0$,  $\sigma_i$, $\eta_{ij}$, $r_{ij} \in  C^1(\mathbb{R},\mathbb{R})$ are nonnegative and Bohr almost periodic, and there exist constants $\dot \sigma_i^+$, $\dot \eta_{ij}^+$, $\dot r_{ij}^+$ such that $\dot \sigma_i \le \dot \sigma_i^+<1$, $\dot \eta_{ij} \le \dot \eta_{ij}^+<1$, $\dot r_{ij} \le \dot r_{ij}^+<1$; moreover, $a_i,  c_{ij}, d_{ij}, U_{ij}, I_i \in B_{PAP}^p(\mathbb{R},\mathbb{R})\cap \mathcal{L}^{\infty}(\mathbb{R},\mathbb{R})$ can be written as
	 $a_i=a_i^1+a_i^0$, $c_{ij}=c_{ij}^1+c_{ij}^0$, $d_{ij}=d_{ij}^1+d_{ij}^0$, $U_{ij}=U_{ij}^1+U_{ij}^0, I_i=I_i^1+I_i^0$, where $a_i^1$, $c_{ij}^1$, $d_{ij}^1$, $U_{ij}^1 \in B_{AP}^p(\mathbb{R},\mathbb{R})\cap \mathcal{L}^{\infty}(\mathbb{R},\mathbb{R})$, $a_i^0$, $c_{ij}^0$, $d_{ij}^0$, $U_{ij}^0,I_i^0 \in \mathcal{O}^p (\mathbb{R},\mathbb{R},\mu,\nu)$.
	\item [$(H_2)$]
	For $j = 1,2, \cdots ,n$, $f_j$, $g_j$, $\delta_j  \in C( \mathbb{R},\mathbb{R} )$,  there exist positive constants $L_j^f, L_j^g$ and $L_j^\delta$ such that for all $x$, $y\in \mathbb{R}$,
	$$| f_j( x ) - f_j( y ) | \le L_j^f | x - y |,\quad
	| g_j( x ) - g_j( y ) | \le L_j^g | x - y |,\quad
	| \delta_j(x) - \delta_j(y) | \le L_j^\delta | x - y |;$$ in addition, $f_j(0) = g_j(0) = \delta_j(0)  = 0$.
	\item [$(H_3)$]
	For $p > 2$, $\frac{1}{p}+\frac{1}{q}=1$, $2^pF_1<1$, and for $p=2$, $4F_2<1$, where
\begin{align*}
  F_1:=&3^{p-1}\max_{1\le i\le n}\bigg\{\bigg(\frac{1}{k_i^-}\bigg)^{p} \bigg[\bigg(\sum_{j = 1}^{n} (c_{ij}^+)^q\bigg)^{\frac{p}{q}} \sum_{j=1}^{n}(L_j^f)^p +  \bigg(\sum_{j = 1}^{n} (d_{ij}^+)^q\bigg)^{\frac{p}{q}} \sum_{j=1}^{n}(L_j^g)^p\bigg]\\
  &+ C_p\bigg(\frac{p-2}{2k_i^-}\bigg)^{\frac{p-2}{2}} \bigg(\sum_{j = 1}^{n} (U_{ij}^+)^q\bigg)^{\frac{p}{q}} \sum_{j=1}^{n}(L_j^\delta)^p\frac{2}{pk_i^-}\bigg\},\\
  F_2:=& 3\max_{1\le i\le n}\bigg\{\bigg(\frac{1}{k_i^-}\bigg)^{2}\bigg[ \bigg(\sum_{j = 1}^{n} (c_{ij}^+)^2\bigg) \sum_{j=1}^{n}(L_j^f)^2+ \bigg(\sum_{j = 1}^{n} (d_{ij}^+)^2\bigg) \sum_{j=1}^{n}(L_j^g)^2\bigg]\\
  &+\frac{1}{2k_i^-}\sum_{j = 1}^{n} (U_{ij}^+)^2 \sum_{j=1}^{n}(L_j^\delta)^2  \bigg\}.
\end{align*}

\item [$(H_4)$]
For $p > 2$, $\frac{1}{p}+\frac{1}{q}=1$, $2^pF'_1<1$, and for $p=2$, $4F'_2<1$, here
\begin{align*}
  F_1':=&3^{p-1}\max_{1\le i\le n}\bigg\{\bigg(\frac{1}{\underline{k}_i^1}\bigg)^{p} \bigg[\bigg(\sum_{j = 1}^{n} (\bar{c}_{ij}^1)^q\bigg)^{\frac{p}{q}} \sum_{j=1}^{n}(L_j^f)^p +  \bigg(\sum_{j = 1}^{n} (\bar{d}_{ij}^1)^q\bigg)^{\frac{p}{q}} \sum_{j=1}^{n}(L_j^g)^p\bigg]\\
  &+ C_p\bigg(\frac{p-2}{2\underline{k}_i^1}\bigg)^{\frac{p-2}{2}} \bigg(\sum_{j = 1}^{n} (\bar{U}_{ij}^1)^q\bigg)^{\frac{p}{q}} \sum_{j=1}^{n}(L_j^\delta)^p\frac{2}{p\underline{k}_i^1}\bigg\},\\
  F_2':=& 3\max_{1\le i\le n}\bigg\{\bigg(\frac{1}{\underline{k}_i^1}\bigg)^{2}\bigg[ \bigg(\sum_{j = 1}^{n} (\bar{c}_{ij}^1)^2\bigg) \sum_{j=1}^{n}(L_j^f)^2+ \bigg(\sum_{j = 1}^{n} (\bar{d}_{ij}^1)^2\bigg) \sum_{j=1}^{n}(L_j^g)^2\bigg]\\
  &+\frac{1}{2\underline{k}_i^1}\sum_{j = 1}^{n} (\bar{U}_{ij}^1)^2 \sum_{j=1}^{n}(L_j^\delta)^2  \bigg\},
\end{align*}
where $k_i^1$, $c_{ij}^1$, $d_{ij}^1$, $U_{ij}^1$ are defined in $(H_1)$.
\item [$(H_5)$]
For $p > 2$, $\frac{1}{p}+\frac{1}{q}=1$,
\begin{align*}
  \Xi_1=&{11}^{p-1}\max_{1\le i\le n}\bigg\{ \bigg(\frac{1}{\underline{k}_i^1}\bigg)^{p}\bigg[\bigg(\sum_{j = 1}^{n} (\bar c_{ij}^1)^q\bigg)^{\frac{p}{q}} \sum_{j=1}^{n}(L_j^f)^p  +{2}^{p-1}   \bigg(\sum_{j = 1}^{n} (\bar d_{ij}^1)^q\bigg)^{\frac{p}{q}} \sum_{j=1}^{n}(L_j^g)^p\\
  &\times \frac{e^{\underline{k}_i^1 \eta_{ij}^+}}{1-\dot \eta_{ij}^+}\bigg]+{2}^{p-1}C_p \bigg(\frac{p-2}{p\underline{k}_i^1}\bigg)^{\frac{p-2}{2}}\frac{2}{p\underline{k}_i^1} \bigg(\sum_{j = 1}^{n} (\bar U_{ij}^1)^q\bigg)^{\frac{p}{q}} \sum_{j=1}^{n}(L_j^\delta)^p\frac{e^{\frac{p}{2}\underline{k}_i^1 r_{ij}^+}}{1-\dot r_{ij}^+} \bigg\}<1;
\end{align*}
for $p=2$,
\begin{align*}
  \Xi_2=&{11} \max_{1\le i\le n}\bigg\{ \bigg(\frac{1}{\underline{k}_i^1}\bigg)^{2}\bigg[\sum_{j = 1}^{n} (\bar c_{ij}^1)^2 \sum_{j=1}^{n}(L_j^f)^2 + 2   \sum_{j = 1}^{n} (\bar d_{ij}^1)^2  \sum_{j=1}^{n}(L_j^g)^2\frac{e^{\underline{k}_i^1 \eta_{ij}^+}}{1-\dot \eta_{ij}^+}\bigg]\\
  &+  \frac{1}{\underline{k}_i^1}\bigg(\sum_{j = 1}^{n} (\bar U_{ij}^1)^2\bigg)  \sum_{j=1}^{n}(L_j^\delta)^2\frac{e^{2\underline{k}_i^1 r_{ij}^+}}{1-\dot r_{ij}^+} \bigg\}<1.
\end{align*}

\item [$(H_6)$]
For $p > 2$, $\frac{1}{p}+\frac{1}{q}=1$,
\begin{align*}
 \beta=&{11}^{p-1} \max_{1\le i\le n}\bigg\{\bigg(\frac{1}{k_i^-}\bigg)^{\frac{p}{q}}\bigg[\bigg(\sum_{j=1}^{n} (\bar c_{ij}^1)^q\bigg)^{\frac{p}{q}} \sum_{j=1}^{n}(L_j^f)^p
+\bigg(\sum_{j=1}^{n} (\bar d_{ij}^1)^q\bigg)^{\frac{p}{q}} \sum_{j=1}^{n}(L_j^g)^p\frac{e^{k_i^-\eta_{ij}^+}}{1-\dot \eta_{ij}^+} \bigg]\\
 &+C_p\bigg(\frac{p-2}{pk_i^-}\bigg)^{\frac{p-2}{2}}\bigg(\sum_{j=1}^{n} (\bar U_{ij}^1)^q \bigg)^{\frac{p}{q}} \sum_{j=1}^{n}(L_j^\delta)^p\frac{e^{\frac{p}{2}k_i^-r_{ij}^+}}{1-\dot r_{ij}^+}\bigg\} <k^-;
\end{align*}
for $p=2$,
\begin{align*}
 \tilde{\beta}=&{11}  \max_{1\le i\le n}\bigg\{\frac{1}{k_i^-}\bigg[ \sum_{j=1}^{n} (\bar c_{ij}^1)^2 \sum_{j=1}^{n}(L_j^f)^2+\sum_{j=1}^{n} (\bar d_{ij}^1)^2 \sum_{j=1}^{n}(L_j^g)^2\frac{e^{k_i^-\eta_{ij}^+}}{1-\dot \eta_{ij}^+} \bigg]\\
 &+ \sum_{j=1}^{n}(\bar U_{ij}^1)^2\sum_{j=1}^{n}(L_j^\delta)^2\frac{e^{2k_i^- r_{ij}^+}}{1-\dot r_{ij}^+}\bigg\} <k^-.
\end{align*}
\end{itemize}

\begin{remark}
Obviously, assumption $(H_5)$ implies assumption $(H_3)$.
\end{remark}

Set
$\mathbb{Y}= \big\{ y\in SB_{BC}^{p}(\mathbb{R},\mathcal{L}^{p}(\Omega,\mathbb{R}^n)) \cap \mathcal{L}^{\infty}(\mathbb{R},\mathcal{L}^{p}(\Omega,\mathbb{R}^n)): y$ is $\mathcal{F}_t$- progressively measurable$\big\}$
equipped with the norm $\|y\|_{\mathbb{Y}}=\mathrm{ess}\sup_{t\in \mathbb{R}}\|y(t)\|_p$. Then it is a Banach space.

Define
$\varphi_0 = \big(\int_{-\infty}^t {e^{-\int_s^t k_1(u)\mathrm{d} u }I_1(s)\mathrm{d} s,\cdots,
\int_{-\infty}^t e^{-\int_s^t k_n (u)\mathrm{d} u } I_n(s)\mathrm{d} s} \big)$, and $\varphi_{0}^{1}$ analogously with $I_i^{1}$.
Take a constant $C \ge \max\{\| \varphi_0 \|_\mathbb{Y},\| \varphi^1_0 \|_\mathbb{Y}\}$. Then, we have
\begin{theorem}\label{th3.1}
	If conditions $(H_1)-(H_3)$ are satisfied, then system \eqref{e1} has a unique  solution in $\mathbb{Y}^*=\{\varphi \in \mathbb{Y}| \|\varphi-\varphi_0\|_\mathbb{Y} \le C\}$.
\end{theorem}
\begin{proof}
	Taking the limit as $t_0 \rightarrow -\infty$, by \eqref{jie}, we have
	\begin{align}\label{equation1}
		x_i(t) =&\int_{-\infty}^t e^{-\int_s^t k_i (u)\mathrm{d} u}\Theta_i(s,x)\mathrm{d} s+\int_{-\infty}^t e^{-\int_s^t k_i (u)\mathrm{d} u}\Xi_{ij} (s,x)\mathrm{d} w_j(s),
	\end{align}
which is also a mild solution of system \eqref{e1}.

	Define an operator $\Phi:\mathbb{Y}^* \to \mathbb{Y}$ by
	$\Phi \varphi =(\Phi_1\varphi,\cdots,\Phi_n\varphi)$, where for $t \in \mathbb{R}$, $i=1,\cdots ,n$,
\begin{align}\label{e43a}
		(\Phi_i\varphi)(t) =&\int_{-\infty}^t e^{-\int_s^t k_i (u)\mathrm{d} u} \Theta_i(s,\varphi) \mathrm{d} s +\int_{-\infty }^t e^{-\int_s^t k_i (u)\mathrm{d} u} \Xi_{ij}(s,\varphi)\mathrm{d} w_j(s).
\end{align}

Firstly, we prove that $\Phi$ maps from $\mathbb{Y}^*$ into $\mathbb{Y}^*$. For any $\varphi \in \mathbb{Y}^*$, note the fact that $\|\varphi\|_\mathbb{Y} \le \|\varphi-\varphi_0\|_\mathbb{Y}+\|\varphi_0\|_\mathbb{Y} \le 2C$. From \eqref{e43a}, we find
\begin{align}\label{3.4}
	\|\Phi \varphi-\varphi_0\|_p^p
	\le &3^{p-1}\max_{1\le i\le n}\bigg\{E\bigg|\int_{-\infty}^t e^{-\int_s^t k_i (u)\mathrm{d} u}\sum_{j = 1}^{n} c_{ij}(s) f_j(\varphi_j(s)) \mathrm{d} s\bigg|^p\bigg\}\nonumber\\
	&+3^{p-1}\max_{1\le i\le n}\bigg\{E\bigg|\int_{-\infty}^t e^{-\int_s^t k_i (u)\mathrm{d} u}\sum_{j = 1}^{n} d_{ij}(s) g_j(\varphi_j(s-\eta_{ij}(s))) \mathrm{d} s\bigg|^p\bigg\}\nonumber\\
	&+3^{p-1}\max_{1\le i\le n}\bigg\{E\bigg|\int_{-\infty}^t e^{-\int_s^t k_i (u)\mathrm{d} u}\sum_{j = 1}^{n} U_{ij}(s) \delta_j(\varphi_j(s-r_{ij}(s))) \mathrm{d} w_j(s)\bigg|^p\bigg\}\nonumber\\
	:=&A_1+A_2+A_3.
\end{align}
For $A_1$, by the H\"older inequality, it holds that
\begin{align}\label{3.5}
	A_1
	\le& 3^{p-1}\max_{1\le i\le n}\bigg\{\bigg[\int_{-\infty}^{t}e^{-\int_{s}^{t}k_i (u)\mathrm{d} u} \mathrm{d} s\bigg]^{\frac{p}{q}}\bigg[\int_{-\infty}^{t}e^{-\int_{s}^{t}k_i (u) \mathrm{d} u}E\bigg|\sum_{j = 1}^{n} c_{ij}(s) f_j(\varphi_j(s))\bigg|^p \mathrm{d} s\bigg]\bigg\}\nonumber\\
	\le &3^{p-1}\max_{1\le i\le n}\bigg\{\bigg(\frac{1}{k_i^-}\bigg)^{\frac{p}{q}} \bigg(\sum_{j = 1}^{n} (c_{ij}^+)^q\bigg)^{\frac{p}{q}} \sum_{j=1}^{n}(L_j^f)^p\int_{-\infty}^{t}e^{-k_i^-(t-s)}E|\varphi_i(s) |^p \mathrm{d} s \bigg\}\nonumber\\
	\le& 3^{p-1}\max_{1\le i\le n}\bigg\{\bigg(\frac{1}{k_i^-}\bigg)^{p} \bigg(\sum_{j = 1}^{n} (c_{ij}^+)^q\bigg)^{\frac{p}{q}} \sum_{j=1}^{n}(L_j^f)^p\bigg\}\|\varphi\|_\mathbb{Y}^p.
\end{align}
Similarly, we can get
\begin{align}
	A_2
	\le 3^{p-1}\max_{1\le i\le n}\bigg\{\bigg(\frac{1}{k_i^-}\bigg)^{p} \bigg(\sum_{j = 1}^{n} (d_{ij}^+)^q\bigg)^{\frac{p}{q}} \sum_{j=1}^{n}(L_j^g)^p\bigg\}\|\varphi\|_\mathbb{Y}^p.
\end{align}
For $p>2$, according to the BDG inequality, we obtain
\begin{align}\label{3.8}
	A_3
	\le &3^{p-1}C_p\max_{1\le i\le n}\bigg\{\bigg[\int_{-\infty}^{t}\Big(e^{-\int_{s}^{t}k_i (u)\mathrm{d} u}\Big)^{ \frac{p}{p-2}} \mathrm{d} s\bigg]^{\frac{p-2}{p}\times \frac{p}{2}}E\bigg[\int_{-\infty}^{t}\Big(e^{-\int_{s}^{t}k_i (u)\mathrm{d} u}\Big)^{ \frac{p}{2}}\nonumber\\
	& \times \bigg(\bigg|\sum_{j = 1}^{n} U_{ij}(s) \delta_j(\varphi_j(s-r_{ij}(s)))\bigg|^2\bigg)^{\frac{p}{2}} \mathrm{d} s\bigg]  \bigg\}\nonumber\\
	\le &3^{p-1}C_p \max_{1\le i\le n}\bigg\{\bigg(\frac{p-2}{pk_i^-}\bigg)^{\frac{p-2}{2}} \bigg(\sum_{j = 1}^{n} (U_{ij}^+)^q\bigg)^{\frac{p}{q}} \sum_{j=1}^{n}(L_j^\delta)^p  \int_{-\infty}^{t}e^{-\frac{p}{2}k_i^-(t-s)}E|\varphi_j(s-r_{ij}(s)) |^p  \mathrm{d} s \bigg\}\nonumber\\
	\le &3^{p-1}C_p\max_{1\le i\le n}\bigg\{\bigg(\frac{p-2}{2k_i^-}\bigg)^{\frac{p-2}{2}} \bigg(\sum_{j = 1}^{n} (U_{ij}^+)^q\bigg)^{\frac{p}{q}} \sum_{j=1}^{n}(L_j^\delta)^p\frac{2}{pk_i^-}\bigg\}\|\varphi \|_\mathbb{Y}^p.
\end{align}
When $p=2$, utilizing the It\^{o} isometry, one can deduce that
\begin{align}\label{3.9}
	A_3' \le &3\max_{1\le i\le n}\bigg\{ E\int_{-\infty}^{t}e^{-2\int_{s}^{t}k_i (u)\mathrm{d} u}\bigg|\sum_{j = 1}^{n} U_{ij}(s) \delta_j(\varphi_j(s-r_{ij}(s)))\bigg|^2 \mathrm{d} s \bigg\}\nonumber\\
	\le &3\max_{1\le i\le n}\bigg\{\sum_{j = 1}^{n} (U_{ij}^+)^2 \sum_{j=1}^{n}(L_j^\delta)^2\frac{1}{2k_i^-}  \bigg\}\|\varphi\|_\mathbb{Y}^2.
\end{align}
Substituting $\eqref{3.5}$-$\eqref{3.9}$ into $\eqref{3.4}$, we have
\begin{align*}
		\|\Phi\varphi-\varphi_0\|_\mathbb{Y}^p
		\le F_1\|\varphi\|_\mathbb{Y}^p \le 2^p F_1C^p \le C^p,\quad (p>2)
    \end{align*}
and
    \begin{align*}
	\|\Phi\varphi -\varphi_0\|_\mathbb{Y}^2
	\le F_2\|\varphi\|_\mathbb{Y}^2 \le4F_2 C^2 \le C^2,\quad (p=2),
\end{align*}
in which $F_1$ and $F_2$ are defined in $(H_3)$. Thus, $\|\Phi\varphi -\varphi_0\|_\mathbb{Y} \le C$ for any $\varphi\in \mathbb{Y}^*$.

Since $\|\Phi\varphi\|_{SB^p}\leq \|\Phi\varphi\|_{\mathbb{Y}}\leq \|\Phi\varphi-\varphi_0\|_\mathbb{Y}+\|\varphi_0\|_\mathbb{Y}\leq 2C$
 for any $\varphi\in \mathbb{Y}^*$, $\Phi\varphi$ is $SB^p$-bounded.

Next, we prove that $(\Phi \varphi)(\cdot)$ is $SB^p$-continuous. For any $\varphi \in \mathbb{Y}^*$ and $\xi>0$, we infer that
\begin{align}\label{q}
&\|(\Phi\varphi)(\cdot+\xi)-\Phi\varphi(\cdot)\|_{SB^p}^p\nonumber\\
\le	&12^{p-1}\max_{1\le i\le n}\bigg\{\overline{\lim_{l\rightarrow \infty}} \frac{1}{2l}\int_{-l}^{l}\bigg[\int_{-\infty}^{t}\Big|e^{-\int_{s}^{t+\xi}k_i (u) \mathrm{d} u}- e^{-\int_s^t k_i (u) \mathrm{d} u}\Big| \mathrm{d} s \bigg]^{\frac{p}{q}}\bigg[\int_{-\infty}^{t}\Big|e^{-\int_{s}^{t+\xi}k_i (u)\mathrm{d} u}\nonumber\\
&- e^{-\int_s^t k_i (u) \mathrm{d} u}\Big|\bigg(E\bigg|\sum_{j = 1}^{n} c_{ij}(s) f_j(\varphi_j(s))\bigg|^p+E\bigg|\sum_{j = 1}^{n} d_{ij}(s) g_j(\varphi_j(s-\eta_{ij}(s)))\bigg|^p\nonumber\\
&+|I_i(s)|^p\bigg)  \mathrm{d} s\bigg] \mathrm{d} t \bigg\}+12^{p-1}\max_{1\le i\le n}\bigg\{\overline{\lim_{l\rightarrow \infty}} \frac{1}{2l}\int_{-l}^{l}\bigg[\int_{t}^{t+\xi}\Big(e^{-\int_{s}^{t+\xi}k_i (u)\mathrm{d} u}\Big)^q \mathrm{d} s \bigg]^{\frac{p}{q}}\nonumber\\
&\times\bigg[\int_{t}^{t+\xi}\bigg(E\bigg|\sum_{j = 1}^{n} c_{ij}(s) f_j(\varphi_j(s))\bigg|^p +E\bigg|\sum_{j = 1}^{n} d_{ij}(s) g_j(\varphi_j(s-\eta_{ij}(s)))\bigg|^p+|I_i(s)|^p\bigg)  \mathrm{d} s\bigg] \mathrm{d} t \bigg\}\nonumber\\
&+4^{p-1}\max_{1\le i\le n}\bigg\{\overline{\lim_{l\rightarrow \infty}} \frac{1}{2l}\int_{-l}^{l}E\bigg| \int_{-\infty}^{t}\Big|e^{-\int_{s}^{t+\xi}k_i (u)\mathrm{d} u}- e^{-\int_s^t k_i (u)\mathrm{d} u}\Big|\sum_{j = 1}^{n} U_{ij}(s)\nonumber\\
&\times  \delta_j(\varphi_j(s-r_{ij}(s))) \mathrm{d} w_j(s)\bigg|^p \mathrm{d} t \bigg\} +4^{p-1}\max_{1\le i\le n}\bigg\{\overline{\lim_{l \rightarrow \infty}} \frac{1}{2l}\int_{-l}^{l}E\bigg| \int_{t}^{t+\xi} e^{-\int_{s}^{t+\xi}k_i (u) \mathrm{d} u}\nonumber\\
&\times  \sum_{j = 1}^{n} U_{ij}(s) \delta_j(\varphi_j(s-r_{ij}(s))) \mathrm{d} w_j(s)\bigg|^p \mathrm{d} t \bigg\} \nonumber\\
	:=&B_1+B_2+B_3+B_4.
\end{align}
According to the H\"older inequality and the Fubini theorem, we can deduce that
\begin{align}\label{3.13}
B_1 \le  &12^{p-1} \max_{1\le i\le n}\bigg\{\overline{\lim_{l\rightarrow \infty}} \frac{1}{2l}\int_{-l}^{l} \bigg(\int_{-\infty}^{t} e^{-k^-_i(t-s)}  \int_{t}^{t+\xi}k_i (u)\mathrm{d} u\mathrm{d} s\bigg)^{\frac{p}{q}}\nonumber\\
&\times \int_{-\infty}^{t}e^{-k_i^-(t-s)} \int_{t}^{t+\xi}k_i (u)\mathrm{d} u\bigg[\bigg(\sum_{j = 1}^{n} (c_{ij}^+)^q\bigg)^{\frac{p}{q}} \sum_{j=1}^{n}(L_j^f)^p E|\varphi_j(s)|^p\nonumber\\
&  +\bigg(\sum_{j = 1}^{n} (d_{ij}^+)^q\bigg)^{\frac{p}{q}}\sum_{j=1}^{n}(L_j^g)^p E|\varphi_j(s-\eta_{ij}(s))|^p +E|I_i(s)|^p\bigg]\mathrm{d} s \mathrm{d} t \bigg\}\nonumber\\
\le &12^{p-1} \max_{1\le i\le n}\bigg\{\bigg(\frac{1}{k_i^-}\bigg)^{p}(k_i^+ \xi)^{p}\bigg[  \bigg(\sum_{j = 1}^{n} (c_{ij}^+)^q\bigg)^{\frac{p}{q}}\sum_{j=1}^{n}(L_j^f)^p \|\varphi\|_\mathbb{Y}^p\nonumber\\
& +\bigg(\sum_{j = 1}^{n} (d_{ij}^+)^q\bigg)^{\frac{p}{q}} \sum_{j=1}^{n}(L_j^g)^p \|\varphi\|_\mathbb{Y}^p +(I_i^+)^p\bigg] \bigg\}
=:M_1\xi^p,\\
B_2\le &12^{p-1} \max_{1\le i\le n}\bigg\{\xi^{\frac{p}{q}} \overline{\lim_{l\rightarrow \infty}} \frac{1}{2l}\int_{-l}^{l} \int_{t}^{t+\xi}\bigg[\bigg(\sum_{j = 1}^{n} (c_{ij}^+)^q\bigg)^{\frac{p}{q}}\sum_{j=1}^{n}(L_j^f)^p \|\varphi\|_\mathbb{Y}^p \nonumber\\
&+\bigg(\sum_{j = 1}^{n} (d_{ij}^+)^q\bigg)^{\frac{p}{q}} \sum_{j=1}^{n}(L_j^g)^p \|\varphi\|_\mathbb{Y}^p +(I^+_i)^p\bigg] \mathrm{d} s \mathrm{d} t \bigg\}\nonumber\\
\le &12^{p-1} \max_{1\le i\le n}\bigg\{\xi^{p} \bigg[\bigg(\sum_{j = 1}^{n} (c_{ij}^+)^q\bigg)^{\frac{p}{q}}\sum_{j=1}^{n}(L_j^f)^p \|\varphi\|_\mathbb{Y}^p \nonumber\\
	&+\bigg(\sum_{j = 1}^{n} (d_{ij}^+)^q\bigg)^{\frac{p}{q}} \sum_{j=1}^{n}(L_j^g)^p \|\varphi\|_\mathbb{Y}^p +(I^+_i)^p\bigg]  \bigg\}
=:M_2\xi^{p}.
\end{align}

For $p>2$, by the BDG inequality and the Fubini theorem, we have
\begin{align}
B_3
	 \le &4^{p-1}C_p \max_{1\le i\le n}\bigg\{\overline{\lim_{l\rightarrow \infty}} \frac{1}{2l}\int_{-l}^{l}\bigg[\int_{-\infty}^{t}e^{- \frac{p}{p-2}k^-_i (t-s) } \bigg(\int_{t}^{t+\xi} k_i (u)\mathrm{d} u\bigg)^{\frac{p}{p-2}} \mathrm{d} s\bigg]^{\frac{p-2}{2}}\nonumber\\
	 &\times \bigg[\int_{-\infty}^{t}e^{-\frac{p}{2}k_i^- (t-s)} \bigg(\int_{t}^{t+\xi} k_i (u)\mathrm{d} u\bigg)^{\frac{p}{2}}\bigg(\sum_{j = 1}^{n} (U_{ij}^+)^q\bigg)^{\frac{p}{q}} \sum_{j=1}^{n}(L_j^\delta)^p\|\varphi\|_\mathbb{Y}^p \mathrm{d} s \bigg]\mathrm{d} t \bigg\}\nonumber\\
	  \le &4^{p-1}C_p \max_{1\le i\le n}\bigg\{\bigg(\frac{p-2}{pk_i^-}\bigg)^{\frac{p-2}{2}}\bigg(\frac{2}{pk_i^-}\bigg)(k_i^+\xi)^p\bigg(\sum_{j = 1}^{n} (U_{ij}^+)^q\bigg)^{\frac{p}{q}} \sum_{j=1}^{n}(L_j^\delta)^p\|\varphi\|_{\mathbb{Y}}^p\bigg\}
=:M_3\xi^p,\\
B_4\leq&4^{p-1}\max_{1\le i\le n}\bigg\{\overline{\lim_{l \rightarrow \infty}} \frac{1}{2l}\int_{-l}^{l}E\bigg[ \int_{t}^{t+\xi} e^{-2\int_{s}^{t+\xi}k_i (u) \mathrm{d} u} \bigg( \sum_{j = 1}^{n} U_{ij}(s) \delta_j(\varphi_j(s-r_{ij}(s)))\bigg)^2 \mathrm{d} s\bigg]^{\frac{p}{2}}\mathrm{d} t \bigg\} \nonumber\\
\le &4^{p-1}C_p \max_{1\le i\le n}\bigg\{\bigg(\sum_{j = 1}^{n} (U_{ij}^+)^q\bigg)^{\frac{p}{q}} \sum_{j=1}^{n}(L_{ij}^\delta)^p\|\varphi\|_{\mathbb{Y}}^p \bigg\}\xi^{\frac{p}{2}}
	:=M_4\xi^{\frac{p}{2}}.
\end{align}
For $p=2$, by the It\^{o} isometry and the Fubini theorem, we deduce that
\begin{align}
B_3' \le& 4 \max_{1\le i\le n}\bigg\{\overline{\lim_{l\rightarrow \infty}} \frac{1}{2l}\int_{-l}^{l}\bigg[\int_{-\infty}^{t}e^{-2 k_i^- (t-s) }\bigg(\int_{t}^{t+\xi}k_i (u)\mathrm{d} u\bigg)^2\nonumber\\
	&\times\sum_{j = 1}^{n} (U_{ij}^+)^2 \sum_{j=1}^{n}(L_j^\delta)^2 E|\varphi_j(s-r_{ij}(s))|^2 \mathrm{d} s \bigg] \mathrm{d} t\bigg\}\nonumber\\
	\le& 4 \max_{1\le i\le n}\bigg\{\bigg(\frac{1}{2k_i^-}\bigg)(k_i^+\xi)^2\sum_{j = 1}^{n} (U_{ij}^+)^2 \sum_{j=1}^{n}(L_j^\delta)^2\|\varphi\|^2_\mathbb{Y} \bigg\}
=:M_3'\xi^2,\\
B_4' \le& 4\max_{1\le i\le n}\bigg\{\sum_{j = 1}^{n} (U_{ij}^+)^2 \sum_{j=1}^{n}(L_j^\delta)^2\|\varphi\|_{\mathbb{Y}}^2\bigg\}\xi
=:M_4'\xi.\label{3.18}
\end{align}

From \eqref{q}-\eqref{3.18}, it follows that
\begin{align*}
	&\|(\Phi\varphi)(\cdot+\xi)-\Phi\varphi(\cdot)\|_{SB^p}^p
	\le M_1\xi^p+M_2\xi^p+M_3\xi^p+M_4\xi^{\frac{p}{2}} \rightarrow 0 \text{ as } \xi \rightarrow 0^+,\,\,(p>2), \\
	&\|(\Phi\varphi)(\cdot+\xi)-\Phi\varphi(\cdot)\|_{SB^p}^2
	\le M_1\xi^2+M_2\xi^{2}+M'_3\xi^2+M'_4\xi \rightarrow 0 \text{ as } \xi \rightarrow 0^+, \,\,(p=2).
\end{align*}
Hence, we conclude that $\Phi$ is $SB^p$-continuous. Therefore, we arrive at $\Phi(\mathbb{Y}^*) \subset \mathbb{Y}^*$.

Next, we prove that $\Phi$ is a contraction mapping. For any $\varphi$, $\psi \in \mathbb{Y}^*$, analogously to the estimation  of \eqref{3.5}-\eqref{3.9}, using  $(H_3)$, for $p>2$, we have
\begin{align*}
	&\|(\Phi\varphi)(t)-(\Phi\psi)(t)\|_p^p \\
	\le
	&3^{p-1}\max_{1\le i\le n}\bigg\{E\bigg|\int_{-\infty}^t e^{-\int_s^t k_i (u)\mathrm{d} u}\sum_{j = 1}^{n} c_{ij}(s) [f_j(\varphi_j(s)) -f_j(\psi_j(s)) ]\mathrm{d} s\bigg|^p\bigg\}\\
	&+3^{p-1}\max_{1\le i\le n}\bigg\{E\bigg|\int_{-\infty}^t e^{-\int_s^t k_i (u)\mathrm{d} u}\sum_{j = 1}^{n} d_{ij}(s) [g_j(\varphi_j(s-\eta_{ij}(s)))-g_j(\psi_j(s-\eta_{ij}(s)))] \mathrm{d} s\bigg|^p\bigg\}\\
	&+3^{p-1}\max_{1\le i\le n}\bigg\{E\bigg|\int_{-\infty}^t e^{-\int_s^t k_i (u)\mathrm{d} u}\sum_{j = 1}^{n} U_{ij}(s) [\delta_j(\varphi_j(s-r_{ij}(s)))-\delta_j(\psi_j(s-r_{ij}(s)))] \mathrm{d} w_j(s)\bigg|^p\bigg\}\\
	\le &F_1\|\varphi-\psi\|_\mathbb{Y}^p,
	\end{align*}
 for $p=2$, one has
 \begin{align*}
 	\|(\Phi\varphi)(t)-(\Phi\psi)(t)\|_\mathbb{Y}^2
 		\le &F_2\|\varphi-\psi \|_\mathbb{Y}^2.
 \end{align*}

 Therefore, by $(H_3)$, for both the cases of $p>2$ and $p =2$, we have
\begin{align*}
	\|(\Phi\varphi)(\cdot)-(\Phi\psi)(\cdot)\|_\mathbb{Y} \le  \frac{1}{2}\|\varphi-\psi\|_\mathbb{Y},
\end{align*}
which implies that $\Phi$ is a contraction mapping. Then, $\Phi$ has a unique fixed point $x$ in $\mathbb{Y}^*(\mathbb{R},\mathcal{L}^p(\Omega, \mathbb{R}^n))$. Thus, \eqref{e1} admits a unique  solution $x \in \mathbb{Y}^*$. The proof is completed.
\end{proof}

Consider the following integral system
\begin{align}\label{e2}
		y_i(t) =& \int_{-\infty}^t e^{-\int_s^t k_i^1 (u) \mathrm{d} u} \bigg[\sum_{j = 1}^{n} c_{ij}^1(s) f_j(y_j(s))+\sum_{j = 1}^{n} d_{ij}^1(s) g_j(y_j(s-\eta_{ij}(s))) \nonumber\\
&+I_i^1(s)\bigg]\mathrm{d} s +\int_{-\infty }^t e^{-\int_s^t k_i^1 (u)\mathrm{d} u} \sum_{j = 1}^{n}U_{ij}^1(s) \delta_j(y_j(s-r_{ij}(s)))\mathrm{d} w_j(s),
	\end{align}
where $i=1,2,\ldots, n$.
Based on the definition of Besicovitch almost periodic stochastic processes and Theorem $\ref{th3.1}$, we prove the following theorem.

\begin{theorem}\label{th3.2}
	If conditions $(H_1)$, $(H_2)$, $(H_4)$ and $(H_5)$ are satisfied, then system \eqref{e2} has a unique $p$-th Besicovitch almost periodic solution in distribution.
\end{theorem}
\begin{proof}
	   According to Theorem \ref{th3.1}, we see that  system $\eqref{e2}$
has a unique solution $y=(y_1(t),\cdots,y_n(t))$ in $\mathbb{Y}^*$. To complete the proof, it remains to show that
 $y$ is $p$-th Besicovitch almost periodic in distribution.
	
 Since $y \in \mathbb{Y}^*$, for any $\varepsilon>0$, there exists $\delta \in (0,\varepsilon)$ such that for any $t_1,t_2\in \mathbb{R}$ with $|t_1-t_2|<\delta$, we have $\|y(\cdot+t_1)-y(\cdot+t_2)\|_{SB^p} < \varepsilon$. By $(H_1)$, there exists a $\tau $ such that
 $$
  | \eta_{ij}(t+\tau)- \eta_{ij}(t) | <\delta,\,
  |r_{ij}(t+\tau)-r_{ij}(t)|<\delta,\,\,
   \| k_i^1(\cdot+\tau)-k_i^1(\cdot) \|_{B^p}<\varepsilon,\, \| I_i^1(\cdot+\tau)-I_i^1(\cdot) \|_{B^p}<\varepsilon,
 $$
$$
 \| c_{ij}^1(\cdot+\tau)- c_{ij}^1(\cdot) \|_{B^p}<\varepsilon,\,\,
 \| d_{ij}^1(\cdot+\tau)-d_{ij}^1(\cdot)\|_{B^p}<\varepsilon,\,\,
\| U_{ij}^1(\cdot+\tau)- U_{ij}^1(\cdot)\|_{B^p} <\varepsilon.
$$
Consequently, we can get
\begin{align}
	&\| y(\cdot-\eta_{ij}(\cdot+\tau))-y(\cdot-\eta_{ij}(\cdot))\|_{B^p}<\varepsilon,\label{3.21}\\
	&\|y(\cdot-r_{ij}(\cdot+\tau))-y(\cdot-r_{ij}(\cdot))\|_{B^p}<\varepsilon.\label{3.22}
\end{align}

By virtue of  \eqref{e2}, one has
\begin{align*}
	&y_i(t+\tau)\\
	=&\int_{-\infty}^t e^{-\int_s^t k_i^1 (u+\tau) \mathrm{d} u} \bigg[\sum_{j = 1}^{n} c_{ij}^1(s+\tau) f_j(y_j(s+\tau))+\sum_{j = 1}^{n} d_{ij}^1(s+\tau) g_j(y_j(s+\tau-\eta_{ij}(s+\tau)))\nonumber\\
& +I_i^1(s+\tau)\bigg]\mathrm{d} s +\int_{-\infty }^t e^{-\int_s^t k_i^1 (u+\tau)\mathrm{d} u} \bigg[\sum_{j = 1}^{n}U_{ij}^1(s+\tau) \delta_j(y_j(s+\tau-r_{ij}(s+\tau)))\bigg] \mathrm{d} \check{w}_j(s),
\end{align*}
  where $i=1,\dots,n$, $\check{w}_j(s)=w_j(s+\tau)-w_j(\tau)$ is a Brownian motion sharing the same distribution with that of $w_j(s)$.

Let us  consider a new process
\begin{align}\label{3.23}
	\hat{y}_i(t)
	=&\int_{-\infty}^t e^{-\int_s^t k_i^1 (u+\tau) \mathrm{d} u} \bigg[\sum_{j = 1}^{n} c_{ij}^1(s+\tau) f_j(\hat{y}_j(s))+\sum_{j = 1}^{n} d_{ij}^1(s+\tau) g_j(\hat{y}_j(s-\eta_{ij}(s+\tau)))\nonumber\\
& +I_i^1(s+\tau)\bigg]\mathrm{d} s +\int_{-\infty }^t e^{-\int_s^t k_i^1 (u+\tau)\mathrm{d} u} \bigg[\sum_{j = 1}^{n}U_{ij}^1(s+\tau) \delta_j(\hat{y}_j(s-r_{ij}(s+\tau)))\bigg] \mathrm{d} w_j(s),
\end{align}
where $i=1,2,\ldots,n.$ Then,
 by the same argument as in the proof of Theorem \ref{th3.1}, one can show that \eqref{3.23} has a unique solution $\hat{y}=(\hat{y}_1,\hat{y}_2,\ldots,\check{y}_n)\in \mathbb{Y}^*$. Clearly, $y_i(t+\tau)$ and $\hat{y}$ have the same distribution.
From $\eqref{e2}$ and \eqref{3.23}, it follows that
\begin{align}\label{3.24}
&\|\hat{y}(\cdot)-y(\cdot)\|_{SB^p}^p  \nonumber\\
\le&  {11}^{p-1} \max_{1\le i\le n}\bigg\{\bigg\|\int_{-\infty}^\cdot e^{-\int_s^\cdot k_i^1 (u+\tau)\mathrm{d} u}\sum_{j=1}^{n}[c_{ij}^1(s+\tau)-c_{ij}^1(s)]f_j(\hat{y}_j(s))\mathrm{d} s\bigg\|_{SB^p}^p \bigg\}\nonumber\\
& +{11}^{p-1} \max_{1\le i\le n}\bigg\{\bigg\|\int_{-\infty}^\cdot e^{-\int_s^\cdot k_i^1 (u+\tau)\mathrm{d} u}\sum_{j=1}^{n}c_{ij}^1(s) [f_j(\hat{y}_j(s))-f_j(y_j(s))]\mathrm{d} s\bigg\|_{SB^p}^p \bigg\}\nonumber\\
	&+ {11}^{p-1} \max_{1\le i\le n}\bigg\{\bigg\|\int_{-\infty}^\cdot \Big[e^{-\int_s^\cdot k_i^1 (u+\tau)\mathrm{d} u}-e^{-\int_s^\cdot k_i^1 (u)\mathrm{d} u}\Big]\sum_{j=1}^{n} c_{ij}^1(s) f_j(y_j(s)) \mathrm{d} s\bigg\|_{SB^p}^p \bigg\}\nonumber\\
    &+{11}^{p-1} \max_{1\le i\le n}\bigg\{\bigg\|\int_{-\infty}^\cdot e^{-\int_s^\cdot k_i^1 (u+\tau)\mathrm{d} u}\sum_{j=1}^{n}[d_{ij}^1(s+\tau)-d_{ij}^1(s)] g_j(\hat{y}_j(s -\eta_{ij}(s+\tau)))\mathrm{d} s\bigg\|_{SB^p}^p \bigg\}\nonumber\\
    &+ {11}^{p-1} \max_{1\le i\le n}\bigg\{\bigg\|\int_{-\infty}^\cdot e^{-\int_s^\cdot k_i^1 (u+\tau)\mathrm{d} u}\sum_{j=1}^{n} d_{ij}^1(s) [g_j(\hat{y}_j(s-\eta_{ij}(s+\tau)))\nonumber\\
    &-g_j(y_j(s- \eta_{ij}(s)))]\mathrm{d} s\bigg\|_{SB^p}^p \bigg\}\nonumber\\	
	&+ {11}^{p-1} \max_{1\le i\le n}\bigg\{\bigg\|\int_{-\infty}^\cdot \Big[e^{-\int_s^\cdot k_i^1 (u+\tau)\mathrm{d} u}-e^{-\int_s^\cdot k_i^1 (u)\mathrm{d} u}\Big]\sum_{j=1}^{n}d_{ij}^1(s) g_j(y_j(s-\eta_{ij}(s)))\mathrm{d} s\bigg\|_{SB^p}^p \bigg\}\nonumber\\
   &+ {11}^{p-1} \max_{1\le i\le n}\bigg\{\bigg\|\int_{-\infty}^\cdot e^{-\int_s^\cdot k_i^1 (u+\tau)\mathrm{d} u} [I_i^1(s+\tau)- I_i^1(s)]\mathrm{d} s\bigg\|_{SB^p}^p \bigg\}\nonumber\\
   &+ {11}^{p-1} \max_{1\le i\le n}\bigg\{\bigg\|\int_{-\infty}^\cdot \Big[e^{-\int_s^\cdot k_i^1 (u+\tau)\mathrm{d} u}-e^{-\int_s^\cdot  k_i^1 (u)\mathrm{d} u}\Big] I_i^1(s)\mathrm{d} s\bigg\|_{SB^p}^p \bigg\}\nonumber\\
   &+ {11}^{p-1} \max_{1\le i\le n}\bigg\{\bigg\|\int_{-\infty}^\cdot e^{-\int_s^\cdot k_i^1 (u+\tau)\mathrm{d} u}\sum_{j=1}^{n}[U_{ij}^1(s+\tau)-U_{ij}^1(s)]\nonumber\\
    &\times\delta_j(\hat{y}_j(s-r_{ij}(s+\tau)))\mathrm{d} w_j(s)\bigg\|_{SB^p}^p \bigg\} +{11}^{p-1} \max_{1\le i\le n}\bigg\{\bigg\|\int_{-\infty}^\cdot e^{-\int_s^\cdot k_i^1 (u+\tau)\mathrm{d} u}\sum_{j=1}^{n}U_{ij}^1(s)\nonumber\\
   &\times [\delta_j(\hat{y}_j(s-r_{ij}(s+\tau)))-\delta_j(y_j(s-r_{ij}(s)))]\mathrm{d} w_j(s)\bigg\|_{SB^p}^p \bigg\}\nonumber\\
   &+ {11}^{p-1} \max_{1\le i\le n}\bigg\{\bigg\|\int_{-\infty}^\cdot \Big[e^{-\int_s^\cdot k_i^1 (u+\tau)\mathrm{d} u}-e^{-\int_s^\cdot k_i^1 (u)\mathrm{d} u}\Big]\nonumber\\
   &\times\sum_{j=1}^{n}U_{ij}^1(s)  \delta_j(y_j(s-r_{ij}(s)))\mathrm{d} w_j(s)\bigg\|_{SB^p}^p \bigg\}
= :\sum_{i=1}^{11}J_i.
\end{align}

According to the H\"older inequality and the Fubini theorem, we can infer that
\begin{align}
J_1
\le &{11}^{p-1}\max_{1\le i\le n}\bigg\{\overline{\lim_{l\rightarrow \infty}} \frac{1}{2l}\int_{-l}^{l}\bigg[\int_{-\infty}^{t}e^{-\int_{s}^{t}k_i^1 (u+\tau)\mathrm{d} u} \mathrm{d} s\bigg]^{\frac{p}{q}}\bigg[\int_{-\infty}^{t}e^{-\int_{s}^{t}k_i^1 (u+\tau) \mathrm{d} u}\nonumber\\
& \times E\bigg|\sum_{j=1}^{n}[c_{ij}^1(s+\tau)-c_{ij}^1(s)]f_j(\hat{y}_j(s))\bigg|^p\mathrm{d} s\bigg]\mathrm{d} t \bigg\}\nonumber\\
	\le &{11}^{p-1}\max_{1\le i\le n}\bigg\{ \bigg(\frac{1}{\underline{k}_i^1}\bigg)^{\frac{p}{q}} \bigg(\sum_{j=1}^{n}(2CL_j^f)^q \bigg)^{\frac{p}{q}}\overline{\lim_{l\rightarrow \infty}} \frac{1}{2l}\int_{0}^{+\infty}\int_{-l}^{l}e^{-\underline{k}_i^1 s}\nonumber\\
	&\times \bigg|\sum_{j=1}^{n}[c_{ij}^1(t-s+\tau)-c_{ij}^1(t-s)]\bigg|^p \mathrm{d} t \mathrm{d} s \bigg\}\nonumber\\
	\le &{11}^{p-1}\max_{1\le i\le n}\bigg\{ \bigg(\frac{1}{\underline{k}_i^1}\bigg)^{\frac{p}{q}} \bigg(\sum_{j=1}^{n}(2CL_j^f)^q \bigg)^{\frac{p}{q}}\overline{\lim_{l\rightarrow \infty}}\int_{0}^{+\infty}\frac{2(l+s)}{2l}\nonumber\\
	&\times
\frac{1}{2(l+s)} \int_{-l-s}^{l+s}e^{-\underline{k}_i^1 s} \bigg|\sum_{j=1}^{n}[c_{ij}^1(t+\tau)-c_{ij}^1(t)]\bigg|^p \mathrm{d} t \mathrm{d} s \bigg\}\nonumber\\
	\le &{11}^{p-1}\max_{1\le i\le n}\bigg\{ \bigg(\frac{1}{\underline{k}_i^1}\bigg)^{\frac{p}{q}} \bigg(\sum_{j=1}^{n}(2CL_j^f)^q \bigg)^{\frac{p}{q}}\int_{0}^{+\infty}e^{-\underline{k}_i^1 s}\mathrm{d} s \nonumber\\
	&\times \overline{\lim_{l\rightarrow \infty}}\frac{1}{2l} \int_{-l}^{l} \bigg|\sum_{j=1}^{n}[c_{ij}^1(t+\tau)-c_{ij}^1(t)]\bigg|^p \mathrm{d} t \bigg\}\nonumber\\
	\le &{11}^{p-1}\max_{1\le i\le n}\bigg\{\bigg(\frac{1}{\underline{k}_i^1}\bigg)^{p}\bigg(\sum_{j=1}^{n}(2CL_j^f)^q \bigg)^{\frac{p}{q}} n^{\frac{p}{q}} \bigg\}\varepsilon^p=:\check{J}_1(p)\varepsilon^p,\nonumber\\
J_2\le &{11}^{p-1}\max_{1\le i\le n}\bigg\{ \bigg(\frac{1}{\underline{k}_i^1}\bigg)^{\frac{p}{q}}\bigg(\sum_{j = 1}^{n} (\bar c_{ij}^1)^q\bigg)^{\frac{p}{q}} \sum_{j=1}^{n}(L_j^f)^p\int_{0}^{+\infty}e^{-\underline{k}_i^1 s }\nonumber\\
	&\times\overline{\lim_{l\rightarrow \infty}} \frac{1}{2l}\int_{-l}^{l}  E|\hat{y}_j(t-s)- y_j(t-s)|^p \mathrm{d} t\mathrm{d} s\bigg\}\nonumber\\
	\le  &{11}^{p-1}\max_{1\le i\le n}\bigg\{ \bigg(\frac{1}{\underline{k}_i^1}\bigg)^{\frac{p}{q}}\bigg(\sum_{j = 1}^{n} (\bar c_{ij}^1)^q\bigg)^{\frac{p}{q}} \sum_{j=1}^{n}(L_j^f)^p\int_{0}^{+\infty}e^{-\underline{k}_i^1 s }\overline{\lim_{l\rightarrow \infty}} \frac{2(l+s)}{2l}\nonumber\\
   &\times \frac{1}{2(l+s)}\int_{-l-s}^{l+s} E|\hat{y}_j(t)- y_j(t)|^p \mathrm{d} t\mathrm{d} s\bigg\}\nonumber\\
	\le &{11}^{p-1}\max_{1\le i\le n}\bigg\{ \bigg(\frac{1}{\underline{k}_i^1}\bigg)^{p}\bigg(\sum_{j = 1}^{n} (\bar c_{ij}^1)^q\bigg)^{\frac{p}{q}} \sum_{j=1}^{n}(L_j^f)^p   \bigg\} \|\hat{y}(\cdot)-y(\cdot)\|_{SB^p}^p.\label{12}
\end{align}

Likewise, one can get
\begin{align*}
	J_4
\le &{11}^{p-1}\max_{1\le i\le n}\bigg\{\bigg(\frac{1}{\underline{k}_i^1}\bigg)^{p}\bigg(\sum_{j=1}^{n}(2CL_j^g)^q \bigg)^{\frac{p}{q}} \bigg\}\varepsilon^p=:\check{J}_4(p)\varepsilon^p,\\
	J_7 \le &{11}^{p-1}\max_{1\le i\le n}\bigg\{\bigg(\frac{1}{\underline{k}_i^1}\bigg)^{p} \bigg\} \varepsilon^p=:\check{J}_7(p)\varepsilon^p.
\end{align*}

Note the fact that
 \begin{align*}
 &\overline{\lim_{l\rightarrow \infty}}\frac{1}{2l} \int_{-l}^{l}\bigg|\int_{-\infty}^{t}e^{-\int_s^t k_i^1 (u+\tau)\mathrm{d} u}-e^{-\int_s^t k_i^1 (u)\mathrm{d} u}\mathrm{d} s\bigg|^p\mathrm{d} t\nonumber\\
 \le &\bigg(\frac{1}{\underline{k}^1_i}\bigg)^{\frac{p}{q}}\overline{\lim_{l\rightarrow \infty}}\frac{1}{2l} \int_{-l}^{l}\int_{-\infty}^{t} e^{-\underline{k}_i^1 (t-s)}\bigg(\int_{s}^{t}|k_i^1 (u+\tau)-k_i^1 (u)|\mathrm{d} u  \bigg)^p\mathrm{d} s\mathrm{d} t\nonumber\\
 \le &\bigg(\frac{1}{\underline{k}^1_i}\bigg)^{\frac{p}{q}}\overline{\lim_{l\rightarrow \infty}}\frac{1}{2l} \int_{-l}^{l}\int_{-\infty}^{t}(t-s)^{p-1} e^{-\underline{k}_i^1 (t-s)}\int_{s}^{t}|k_i^1 (u+\tau)-k_i^1 (u)|^p\mathrm{d} u \mathrm{d} s\mathrm{d} t\nonumber\\
 \le &\bigg(\frac{1}{\underline{k}^1_i}\bigg)^{\frac{p}{q}}\int_{0}^{\infty}s^{p-1} e^{-\underline{k}_i^1s}\overline{\lim_{l\rightarrow \infty}}\frac{1}{2l} \int_{-l}^{l}\int_{t-s}^{t}|k_i^1 (u+\tau)-k_i^1 (u)|^p\mathrm{d} u \mathrm{d} t\mathrm{d} s\nonumber\\
   \le &\bigg(\frac{1}{\underline{k}^1_i}\bigg)^{\frac{p}{q}}\int_{0}^{\infty}s^{p-1} e^{-\underline{k}_i^1s}\int_{0}^{s}\overline{\lim_{l\rightarrow \infty}}\frac{1}{2l} \int_{-l}^{l}|k_i^1 (u+t-s+\tau)-k_i^1 (u+t-s)|^p \mathrm{d} t\mathrm{d} u\mathrm{d} s\nonumber\\
   \le&\bigg(\frac{1}{\underline{k}^1_i}\bigg)^{\frac{p}{q}}\frac{\Gamma(p+1)}{(\underline{k}_i^1)^{p+1}}\varepsilon^p.
 \end{align*}
We can obtain that
\begin{align*}
	J_3\leq&{11}^{p-1} \max_{1\le i\le n}\bigg\{\overline{\lim_{l\rightarrow \infty}}\frac{1}{2l} \int_{-l}^{l}\bigg(\int_{-\infty}^t \Big|e^{-\int_s^t k_i^1 (u+\tau)\mathrm{d} u}-e^{-\int_s^t k_i^1 (u)\mathrm{d} u}\Big|\mathrm{d}s\bigg)^{\frac{p}{q}}\nonumber\\
   &\times \int_{-\infty}^t \Big|e^{-\int_s^t k_i^1 (u+\tau)\mathrm{d} u}-e^{-\int_s^t k_i^1 (u)\mathrm{d} u}\Big|\bigg(\sum_{j=1}^{n} |c_{ij}^1(s) f_j(y_j(s))|\bigg)^p \mathrm{d} s  \mathrm{d}t\bigg\}\nonumber\\
 \le&{11}^{p-1} \max_{1\le i\le n}\bigg\{\bigg(\sum_{j = 1}^{n} (\bar c_{ij}^1)^q\bigg)^{\frac{p}{q}} \sum_{j=1}^{n}(2CL_j^f)^p \bigg(\frac{1}{\underline{k}^1_i}\bigg)^{\frac{p}{q}}\frac{\Gamma(p+1)}{(\underline{k}_i^1)^{p+1}}\bigg\}\varepsilon^p=:\check{J}_3(p)\varepsilon^p,\\
	J_6 \le&{11}^{p-1} \max_{1\le i\le n}\bigg\{\bigg(\sum_{j = 1}^{n} (\bar d_{ij}^1)^q\bigg)^{\frac{p}{q}} \sum_{j=1}^{n}(2CL_j^g)^p \bigg(\frac{1}{\underline{k}^1_i}\bigg)^{\frac{p}{q}}\frac{\Gamma(p+1)}{(\underline{k}_i^1)^{p+1}}\bigg\}\varepsilon^p=:\check{J}_6(p)\varepsilon^p,\\
	J_8 \le&{11}^{p-1} \max_{1\le i\le n}\bigg\{(\bar I_i^1)^p \bigg(\frac{1}{\underline{k}^1_i}\bigg)^{\frac{p}{q}}\frac{\Gamma(p+1)}{(\underline{k}_i^1)^{p+1}}  \bigg\}\varepsilon^p=:\check{J}_8(p)\varepsilon^p.
\end{align*}

As for $J_5$, we derive that
\begin{align*}
	J_5
\le &{22}^{p-1} \max_{1\le i\le n}\bigg\{\bigg\|\int_{-\infty}^\cdot e^{-\int_s^\cdot k_i^1 (u+\tau) \mathrm{d} u}\sum_{j=1}^{n}d_{ij}^1(s) [g_j(\hat{y}_j(s-\eta_{ij}(s+\tau))) \\
	&-g_j( y_j(s-\eta_{ij}(s+\tau)))]\mathrm{d} s\bigg\|_{SB^p}^p \bigg\}+{22}^{p-1} \max_{1\le i\le n}\bigg\{\bigg\|\int_{-\infty}^\cdot e^{-\int_s^\cdot k_i^1 (u+\tau)\mathrm{d} u}\sum_{j=1}^{n}d_{ij}^1(s)\nonumber\\
	&\times [g_j(y_j(s-\eta_{ij}(s+\tau)))-g_j(y_j(s-\eta_{ij}(s)))]\mathrm{d} s\bigg\|_{SB^p}^p \bigg\}
=:D_1+D_2.
\end{align*}
In the same way as the method used to estimate $\eqref{12}$, we have
\begin{align*}
	D_1
\le&{22}^{p-1} \max_{1\le i\le n}\bigg\{\bigg(\frac{1}{\underline{k}_i^1}\bigg)^{\frac{p}{q}} \bigg(\sum_{j = 1}^{n} (\bar d_{ij}^1)^q\bigg)^{\frac{p}{q}} \sum_{j=1}^{n}(L_j^g)^p\overline{\lim_{l\rightarrow \infty}} \frac{1}{2l}\int_{-l}^{l}\int_{-\infty}^t e^{-  \underline{k}_i^1 (t-s)}\nonumber\\
&\times \|\hat{y}_j(s-\eta_{ij}(s+\tau))- y_j(s-\eta_{ij}(s+\tau))\|^p\mathrm{d} s \mathrm{d}t \bigg\}\nonumber\\
 \le &{22}^{p-1} \max_{1\le i\le n}\bigg\{\bigg(\frac{1}{\underline{k}_i^1}\bigg)^{p}  \bigg(\sum_{j = 1}^{n} (\bar d_{ij}^1)^q\bigg)^{\frac{p}{q}} \sum_{j=1}^{n}(L_j^g)^p\frac{e^{\underline{k}_i^1 \eta_{ij}^+}}{1-\dot \eta_{ij}^+}\bigg\}\|\hat{y}(\cdot)- y(\cdot)\|_{SB^p}^p.
\end{align*}
Based on \eqref{3.21}, one can get
\begin{align*}
	D_2 \le &{22}^{p-1} \max_{1\le i\le n}\bigg\{\bigg(\frac{1}{\underline{k}_i^1}\bigg)^{p}  \bigg(\sum_{j = 1}^{n} (\bar d_{ij}^1)^q\bigg)^{\frac{p}{q}} \sum_{j=1}^{n}(L_j^g)^p \bigg\}\varepsilon^p=:\check{J}_5(p)\varepsilon^p.
\end{align*}

When $p>2$, as to $J_9$,  using the BDG inequality and the H\"{o}lder inequality, one can get that
\begin{align}\label{3.33}
J_9 \le &{11}^{p-1}  C_p\max_{1\le i\le n}\bigg\{\bigg(\frac{p-2}{p\underline{k}_i^1}\bigg)^{\frac{p-2}{2}}\bigg(\sum_{j=1}^{n}(2CL_j^\delta)^q \bigg)^{\frac{p}{q}}\frac{2}{p\underline{k}_i^1} n^{\frac{p}{q}} \bigg\}\varepsilon^p=:\check{J}_9(p)\varepsilon^p.
\end{align}

For $J_{10}$, we have
\begin{align}
	J_{10}\le &{22}^{p-1} \max_{1\le i\le n}\bigg\{\overline{\lim_{l\rightarrow \infty}}\frac{1}{2l} \int_{-l}^{l}E\bigg|\int_{-\infty}^t e^{-\int_s^t k_i^1 (u+\tau) \mathrm{d} u}\sum_{j=1}^{n}U_{ij}^1(s)\nonumber\\
	&\times [\delta_j(\hat{y}_j(s-r_{ij}(s+\tau)))-\delta_j(y_j(s- r_{ij}(s+\tau)))]\mathrm{d} w_j(s)\bigg|^p \mathrm{d} t\bigg\}\nonumber\\
	&+{22}^{p-1} \max_{1\le i\le n}\bigg\{\overline{\lim_{l\rightarrow \infty}}\frac{1}{2l} \int_{-l}^{l}E\bigg|\int_{-\infty}^t e^{-\int_s^t k_i^1 (u+\tau)\mathrm{d} u}\sum_{j=1}^{n}U_{ij}^1(s)\nonumber\\
	&\times [\delta_j(y_j(s-r_{ij}(s+\tau)))-\delta_j(y_j(s-r_{ij}(s)))]\mathrm{d} w_j(s)\bigg|^p \mathrm{d} t\bigg\}
	:=D_3+D_4.\nonumber
\end{align}
When $p>2$,  by the BDG inequality, the H\"older inequality and the Fubini theorem, we can get
\begin{align*}
	D_3 \le &{22}^{p-1}C_p \max_{1\le i\le n}\bigg\{\bigg(\frac{p-2}{p\underline{k}_i^1}\bigg)^{\frac{p-2}{2}}\frac{2}{p\underline{k}_i^1} \bigg(\sum_{j = 1}^{n} (\bar U_{ij}^1)^q\bigg)^{\frac{p}{q}} \sum_{j=1}^{n}(L_j^\delta)^p\frac{e^{\frac{p}{2}\underline{k}_i^1 r_{ij}^+}}{1-\dot r_{ij}^+} \bigg\}\|\hat{y}(t)-y(t)\|_{SB^p}^p.
\end{align*}
Similar to \eqref{3.33}, by \eqref{3.22}, we obtain that
\begin{align*}
	D_4 \le &{22}^{p-1}C_p \max_{1\le i\le n}\bigg\{\bigg(\frac{p-2}{p\underline{k}_i^1}\bigg)^{\frac{p-2}{2}} \bigg(\sum_{j = 1}^{n} (\bar U_{ij}^1)^q\bigg)^{\frac{p}{q}} \sum_{j=1}^{n}(L_j^\delta)^p\frac{2}{p\underline{k}_i^1} \bigg\}\varepsilon^p=:\check{J}_{10}(p)\varepsilon^p.
\end{align*}

When $p>2$, in view of the BDG inequality and the H\"{o}lder inequality, we can deduce that
\begin{align*}
	J_{11}
	\le&{11}^{p-1}C_p \max_{1\le i\le n}\bigg\{\bigg(\sum_{j = 1}^{n} (\bar U_{ij}^1)^q\bigg)^{\frac{p}{q}} \sum_{j=1}^{n}(2CL_j^\delta)^p\overline{\lim_{l\rightarrow \infty}}\frac{1}{2l} \int_{-l}^{l}\bigg[\int_{-\infty}^t e^{-2\underline{k}_i^1(t-s)} \nonumber\\
	&\times \bigg(\int_{s}^{t}|k_i^1 (u+\tau)-k_i^1 (u)| \mathrm{d} u\bigg)^2 \mathrm{d} s \bigg]^{\frac{p}{2}} \mathrm{d} t\bigg\} \nonumber\\
	\le&{11}^{p-1}C_p \max_{1\le i\le n}\bigg\{\bigg(\sum_{j = 1}^{n} (\bar U_{ij}^1)^q\bigg)^{\frac{p}{q}} \sum_{j=1}^{n}(2CL_j^\delta)^p \bigg(\frac{p-2}{\underline{k}_i^1p}\bigg)^{\frac{p-2}{2}}\overline{\lim_{l\rightarrow \infty}}\frac{1}{2l} \int_{-l}^{l}\int_{-\infty}^t e^{-\frac{\underline{k}_i^1p}{2}(t-s)} \nonumber\\
	&\times
	 \bigg(\int_{s}^{t}|k_i^1 (u+\tau)-k_i^1 (u)| \mathrm{d} u\bigg)^p \mathrm{d} s \mathrm{d} t\bigg\} \nonumber\\
	 \le&{11}^{p-1}C_p \max_{1\le i\le n}\bigg\{\bigg(\sum_{j = 1}^{n} (\bar U_{ij}^1)^q\bigg)^{\frac{p}{q}} \sum_{j=1}^{n}(2CL_j^\delta)^p \bigg(\frac{p-2}{\underline{k}_i^1p}\bigg)^{\frac{p-2}{2}} \nonumber\\
	 &\times\overline{\lim_{l\rightarrow \infty}}\frac{1}{2l} \int_{-l}^{l}\int_{-\infty}^t (t-s)^{p-1} e^{-\frac{\underline{k}_i^1}{2}(t-s)}
	 \int_{s}^{t}|k_i^1 (u+\tau)-k_i^1 (u)|^p \mathrm{d} u \mathrm{d} s \mathrm{d} t\bigg\} \nonumber\\
\le&{11}^{p-1}C_p \max_{1\le i\le n}\bigg\{\bigg(\sum_{j = 1}^{n} (\bar U_{ij}^1)^q\bigg)^{\frac{p}{q}} \sum_{j=1}^{n}(2CL_j^\delta)^p \bigg(\frac{p-2}{\underline{k}_i^1p}\bigg)^{\frac{p-2}{2}}  \nonumber\\
	 &\times\int_{0}^\infty s^{p-1} e^{-\frac{\underline{k}_i^1}{2} s}\int_{0}^{s}\overline{\lim_{l\rightarrow \infty}}
\frac{1}{2l} \int_{-l}^{l}|k_i^1 (u+t-s+\tau)-k_i^1 (u+t-s)|^p \mathrm{d} t\mathrm{d} u\mathrm{d} s\bigg\} \nonumber\\
	 \le&{11}^{p-1}C_p \max_{1\le i\le n}\bigg\{\bigg(\sum_{j = 1}^{n} (\bar U_{ij}^1)^q\bigg)^{\frac{p}{q}} \sum_{j=1}^{n}(2CL_j^\delta)^p \bigg(\frac{p-2}{\underline{k}_i^1p}\bigg)^{\frac{p-2}{2}} \frac{2^{p+1}\Gamma(p+1)}{(\underline{k}_i^1)^{p+1}}   \bigg\}\varepsilon^p=:\check{J}_{11}(p)\varepsilon^p.
\end{align*}

When $p=2$, according to the It\^{o} isometry, we obtain
\begin{align*}
	J'_9 \le& 11\max_{1 \le i \le n}\bigg\{\frac{1}{2\underline{k}_i^1}\sum_{j=1}^{n}(2CL_j^\delta)^2n \bigg\}\varepsilon^2=:\check{J}'_9(2)\varepsilon^2,\\
	J'_{10} \le &22 \max_{1\le i\le n}\bigg\{\overline{\lim_{l\rightarrow \infty}}\frac{1}{2l} \int_{-l}^{l}E\bigg|\int_{-\infty}^t e^{-\int_s^t k_i^1 (u+\tau) \mathrm{d} u}\sum_{j=1}^{n}U_{ij}^1(s)\nonumber\\
	&\times [\delta_j(\hat{y}_j(s-r_{ij}(s+\tau)))-\delta_j(y_j(s- r_{ij}(s+\tau)))]\mathrm{d} w_j(s)\bigg|^2 \mathrm{d} t\bigg\}\nonumber\\
	&+22 \max_{1\le i\le n}\bigg\{\overline{\lim_{l\rightarrow \infty}}\frac{1}{2l} \int_{-l}^{l}E\bigg|\int_{-\infty}^t e^{-\int_s^t k_i^1 (u+\tau)\mathrm{d} u}\sum_{j=1}^{n}U_{ij}^1(s)\nonumber\\
	&\times [\delta_j(y_j(s-r_{ij}(s+\tau)))-\delta_j(y_j(s-r_{ij}(s)))]\mathrm{d} w_j(s)\bigg|^2 \mathrm{d} t\bigg\}
	:=D'_3+D'_4.
\end{align*}

For $D'_3$, one can get
\begin{align*}
	D'_3 \le &22\max_{1\le i\le n}\bigg\{ \frac{1}{2\underline{k}_i^1}\bigg(\sum_{j = 1}^{n} (\bar U_{ij}^1)^2\bigg)  \sum_{j=1}^{n}(L_j^\delta)^2\frac{e^{2\underline{k}_i^1 r_{ij}^+}}{1-\dot r_{ij}^+} \bigg\}\|\hat{y}(t)-y(t)\|_{SB^p}^2.
\end{align*}

Based on \eqref{3.22}, we get that
\begin{align*}
	D'_4 \le &22 \max_{1\le i\le n}\bigg\{ \frac{1}{2\underline{k}_i^1}\sum_{j = 1}^{n} (\bar U_{ij}^1)^2 \sum_{j=1}^{n}(L_j^\delta)^2\bigg\}\varepsilon^2=:\check{J}'_{10}(2)\varepsilon^2.
\end{align*}

By \eqref{3.22},  the It\^{o} isometry and the Cauchy--Schwarz inequality, we infer that
\begin{align*}
J'_{11}\le&11 \max_{1\le i\le n}\bigg\{\sum_{j = 1}^{n} (\bar U_{ij}^1)^2 \sum_{j=1}^{n}(2CL_j^\delta)^2\overline{\lim_{l\rightarrow \infty}}\frac{1}{2l} \int_{-l}^{l}\int_{-\infty}^t (t-s) e^{-2\underline{k}_i^1(t-s)} \nonumber\\
	&\times  \int_{s}^{t}|k_i^1 (u+\tau)-k_i^1 (u)|^2 \mathrm{d} u \mathrm{d} s \mathrm{d} t\bigg\} \nonumber\\
\le&11 \max_{1\le i\le n}\bigg\{\sum_{j = 1}^{n} (\bar U_{ij}^1)^2 \sum_{j=1}^{n}(2CL_j^\delta)^2\int_{0}^\infty s e^{-2\underline{k}_i^1 s}   \int_{0}^{s}\overline{\lim_{l\rightarrow \infty}}\frac{1}{2l}\nonumber\\
 &\times\int_{-l}^{l}|k_i^1 (u+t-s+\tau)-k_i^1 (u+t-s)|^2 \mathrm{d} t\mathrm{d} u \mathrm{d} s\bigg\} \nonumber\\
 \le&11 \max_{1\le i\le n}\bigg\{\sum_{j = 1}^{n} (\bar U_{ij}^1)^2 \sum_{j=1}^{n}(2CL_j^\delta)^2 \frac{2}{(2\underline{k}_i^1)^3}\bigg\}\varepsilon^2=:\check{J}'_{11}(2)\varepsilon^2.
\end{align*}

Therefore, from \eqref{3.24} and the estimates of $J_i$s, for $p>2$, it follows that
\begin{align*}
	\overline{\lim_{l\rightarrow \infty}}\frac{1}{2l} \int_{-l}^{l}E|\hat{y}(t)-y(t)|_0^p \mathrm{d} t
	\le \Lambda_1\varepsilon^p+\Xi_1 \overline{\lim_{l\rightarrow \infty}}\frac{1}{2l} \int_{-l}^{l}E|\hat{y}(t)-y(t)|_0^p \mathrm{d} t,
\end{align*}
where $\Lambda_1:=\sum_{k\in\{1,3,4,\ldots,11\}}\check{J}_k(p)$ and $\Xi_1$ is defined in $(H_5)$.
 By $(H_5)$, we know that $\Xi_1<1$. Hence, we have
\begin{align}\label{ccl1}
	\overline{\lim_{l\rightarrow \infty}}\frac{1}{2l} \int_{-l}^{l}E|\hat{y}(t)- y(t)|_0^p \mathrm{d} t \le \frac{\Lambda_1\varepsilon^p}{1-\Xi_1}.
\end{align}

Likewise, for $p=2$, one can get
\begin{align*}
	\overline{\lim_{l\rightarrow \infty}}\frac{1}{2l} \int_{-l}^{l}E|\hat{y}(t)-y(t)|_0^2 \mathrm{d} t
	\le \Lambda_2\varepsilon^2+\Xi_2\overline{\lim_{l\rightarrow \infty}}\frac{1}{2l} \int_{-l}^{l}E|\hat{y}(t)-y(t)|_0^2 \mathrm{d} t,
\end{align*}
where $\Lambda_2=:\check{J}_1(2) +\check{J}_1(3)+\check{J}_4(2)+\check{J}_5(2)+\check{J}_6(2)+\check{J}_7(2)+\check{J}_8(2)+\check{J}'_9(2)+\check{J}'_{10}(2)+\check{J}'_{11}(2)$ and $\Xi_2$ is mentioned in $(H_5)$.

 Based on $(H_5)$, we have $\Xi_2<1$.
Consequently,  we can infer that
\begin{align}\label{ccl2}
	\overline{\lim_{l\rightarrow \infty}}\frac{1}{2l} \int_{-l}^{l}E|\hat{y}(t)- y(t)|_0^2 \mathrm{d} t < \frac{\Lambda_2\varepsilon^2}{1-\Xi_2}.
\end{align}

Since $y_i(t+\tau)(t)$ and $\hat{y}_i(t)$ have the same distribution for all $t\in \mathbb{R}$, we deduce that
\begin{align*}
  \overline{\lim_{l \rightarrow \infty}} \frac{1}{2l}\int_{-l}^{l} d_{BL}^p(P \circ [y_i(t+\tau)]^{-1}, P \circ [y_i(t)]^{-1}) \mathrm{d} t
  =&\overline{\lim_{l \rightarrow \infty}} \frac{1}{2l}\int_{-l}^{l} d_{BL}^p(P \circ [\hat{y}_i(t)]^{-1}, P \circ [y_i(t)]^{-1}) \mathrm{d} t \\
   \leq&\overline{\lim_{l \rightarrow \infty}} \frac{1}{2l}\int_{-l}^{l}\bigg(  \int_\Omega |\hat{y}_i(t) - y_i(t)|\mathrm{d} P\bigg|\bigg)^p \mathrm{d}t\\
   \leq&\overline{\lim_{l \rightarrow \infty}} \frac{1}{2l}\int_{-l}^{l}  E |\hat{y}_i(t) - y_i(t)|^p \mathrm{d}t,\quad (p\geq 2),
\end{align*}
which combined with \eqref{ccl1} and \eqref{ccl2}, and Definition \ref{bp1} implies that $y(t)$ is $p$-th Besicovitch almost periodic solution in distribution.
 The proof is finished.
\end{proof}

According to the definition of Besicovitch  pseudo almost periodic stochastic processes, Theorems $\ref{th3.1}$ and $\ref{th3.2}$,  we prove the following result.

\begin{theorem}Let $\mu,\nu\in \mathcal{M}_p$, $(A_1)$, $(A_2)$ hold.
	If conditions $(H_1)-(H_4)$ and $(H_6)$ are satisfied, then system \eqref{e1} has a unique  $p$-th Besicovitch $(\mu,\nu )$-pseudo almost periodic solution in distribution.
\end{theorem}
\begin{proof}
		By Theorem $\ref{th3.1}$, the system $\eqref{e1}$ admits a unique solution $x\in \mathbb{Y}^*$ satisfying \eqref{equation1}, and by Theorem \ref{th3.2}, we obtain that   system $\eqref{e2}$ has a unique Besicovitch almost periodic solution in distribution $y\in \mathbb{Y}^*$.
		Now, let $x-y=z$. Then, $x=y+z$. Hence, to complete the proof, we only need to prove that $z \in\mathcal{O}_D^p(\mathbb{R},\mathcal{L}^p(\Omega,\mathbb{R}^n),\mu,\nu)$.
	
From \eqref{equation1} and \eqref{e2}, we have
\begin{align}\label{3.58}
& E|z(t)|_0^p \nonumber\\
\le &{11}^{p-1}\max_{1\le i\le n}\bigg\{E\bigg|\int_{-\infty}^t e^{-\int_s^t k_i (u) \mathrm{d} u}\sum_{j=1}^{n}[c_{ij}(s)-c_{ij}^1(s)] f_j(x_j(s)) \mathrm{d} s\bigg|^p \bigg\}\nonumber\\
&+{11}^{p-1}\max_{1\le i\le n}\bigg\{E\bigg|\int_{-\infty}^t e^{-\int_s^t k_i (u) \mathrm{d} u}\sum_{j=1}^{n}c_{ij}^1(s) [f_j(x_j(s))- f_j(y_j(s))] \mathrm{d} s\bigg|^p \bigg\}\nonumber\\
&+{11}^{p-1}\max_{1\le i\le n}\bigg\{E\bigg|\int_{-\infty}^t \Big[e^{-\int_s^t k_i (u) \mathrm{d} u}-e^{-\int_s^t k_i^1 (u) \mathrm{d} u}\Big]\sum_{j=1}^{n}c_{ij}^1(s) f_j(y_j(s))\mathrm{d} s\bigg|^p \bigg\}\nonumber\\
&+{11}^{p-1}\max_{1\le i\le n}\bigg\{E\bigg|\int_{-\infty}^t e^{-\int_s^t k_i (u) \mathrm{d} u}\sum_{j=1}^{n}[d_{ij}(s)-d_{ij}^1(s)]g_j(x_j(s-\eta_{ij}(s))) \mathrm{d} s\bigg|^p\bigg\}\nonumber\\
&+{11}^{p-1}\max_{1\le i\le n}\bigg\{E\bigg|\int_{-\infty}^t e^{-\int_s^t k_i (u) \mathrm{d} u}\sum_{j=1}^{n}d_{ij}^1(s) [g_j(x_j(s-\eta_{ij}(s)))- g_j(y_j(s-\eta_{ij}(s)))]\mathrm{d} s\bigg|^p\bigg\}\nonumber\\
&+{11}^{p-1}\max_{1\le i\le n}\bigg\{E\bigg|\int_{-\infty}^t \Big[e^{-\int_s^t k_i (u) \mathrm{d} u}-e^{-\int_s^t k_i^1 (u) \mathrm{d} u}\Big]\sum_{j=1}^{n}d_{ij}^1(s)g_j(y_j(s-\eta_{ij}(s))) \mathrm{d} s\bigg|^p\bigg\}\nonumber\\
&+{11}^{p-1}\max_{1\le i\le n}\bigg\{E\bigg|\int_{-\infty}^t e^{-\int_s^t k_i (u)\mathrm{d} u}[I_i(s)-I_i^1(s)] \mathrm{d} s\bigg|^p\bigg\}\nonumber\\
&+{11}^{p-1}\max_{1\le i\le n}\bigg\{E\bigg|\int_{-\infty}^t \Big[e^{-\int_s^t k_i (u)\mathrm{d} u}-e^{-\int_s^t k_i^1 (u)\mathrm{d} u}\Big] I_i^1(s) \mathrm{d} s\bigg|^p\bigg\}\nonumber\\
&+{11}^{p-1}\max_{1\le i\le n}\bigg\{E\bigg|\int_{-\infty}^t e^{-\int_s^t k_i (u) \mathrm{d} u} \sum_{j=1}^{n}[U_{ij}(s)-U_{ij}^1(s)] \delta_j(x_j(s-r_{ij}(s))) \mathrm{d} w_j(s)\bigg|^p\bigg\}\nonumber\\
 &+{11}^{p-1}\max_{1\le i\le n}\bigg\{E\bigg|\int_{-\infty}^t e^{-\int_s^t k_i (u)\mathrm{d} u} \sum_{j=1}^{n}U_{ij}^1(s)[\delta_j(x_j(s-r_{ij}(s))) \nonumber\\
&- \delta_j(y_j(s-r_{ij}(s)))]\mathrm{d} w_j(s)\bigg|^p\bigg\}+{11}^{p-1}\max_{1\le i\le n}\bigg\{E\bigg|\int_{-\infty}^t \Big[e^{-\int_s^t k_i (u)\mathrm{d} u}-e^{-\int_s^t k_i^1 (u)\mathrm{d} u}\Big] \nonumber\\
&\times \sum_{j=1}^{n}U_{ij}^1(s)\delta_j(y_j(s-r_{ij}(s)))\mathrm{d} w_j(s)\bigg|^p\bigg\}
:=\sum_{i=1}^{11}G_i.
	\end{align}

By the H\"older inequality, we get
 \begin{align}
 G_1  \le &{11}^{p-1}\max_{1\le i\le n} \bigg\{\bigg[\int_{-\infty}^{t}e^{-\int_{s}^{t}k_i(u)\mathrm{d} u} \mathrm{d} s\bigg]^{\frac{p}{q}}\bigg[\int_{-\infty}^{t}e^{-\int_{s}^{t}k_i(u)}E\bigg|\sum_{j=1}^{n}c_{ij}^0(s)f_j(x_j(s))\bigg|^p\mathrm{d} s \bigg]\bigg\}\nonumber\\
 \le &{11}^{p-1}\max_{1\le i\le n}\bigg\{\bigg(\frac{1}{k_i^-}\bigg)^{\frac{p}{q}}\bigg(\sum_{j=1}^{n}(2CL_j^f)^q\bigg)^{\frac{p}{q}} \int_{-\infty}^t e^{-k_i^-(t-s)} \sum_{j=1}^{n} |c_{ij}^0(s)|^p \mathrm{d} s\bigg\}:=\alpha_1(t,p),\\
G_2 \le &{11}^{p-1} \max_{1\le i\le n}\bigg\{\bigg(\frac{1}{k_i^-}\bigg)^{\frac{p}{q}}\bigg(\sum_{j=1}^{n} (\bar c_{ij}^1)^q\bigg)^{\frac{p}{q}} \sum_{j=1}^{n}(L_j^f)^p
 \int_{-\infty}^t e^{-k_i^-(t-s)} E|z(s)|_0^p \mathrm{d} s\bigg\},\\
 G_3\le &{11}^{p-1} \max_{1\le i\le n}\bigg\{E\bigg(\int_{-\infty}^{t} e^{-\underline{k}^1_i (t-s)}\int_{s}^{t}|k_i^0 (u)|\mathrm{d} u \bigg|\sum_{j=1}^{n} c_{ij}^1(s)f_j(y_j(s))\mathrm{d} s\bigg|\bigg)^p\bigg\}\nonumber\\
 \le&{11}^{p-1} \max_{1\le i\le n}\bigg\{\bigg(\frac{1}{\underline{k}_i^1}\bigg)^{\frac{p}{q}}\bigg(\sum_{j=1}^{n} (\bar c_{ij}^1)^q\bigg)^{\frac{p}{q}} \sum_{j=1}^{n}(2CL_j^f)^p\int_{-\infty}^{t} e^{-\underline{k}^1_i (t-s)} \bigg(\int_{s}^{t}|k_i^0 (u)|\mathrm{d} u \bigg)^p\mathrm{d} s\bigg\}\nonumber\\
 \le& {11}^{p-1} \max_{1\le i\le n}\bigg\{\bigg(\frac{1}{\underline{k}_i^1}\bigg)^{\frac{p}{q}}\bigg(\sum_{j=1}^{n} (\bar c_{ij}^1)^q\bigg)^{\frac{p}{q}} \sum_{j=1}^{n}(2CL_j^f)^p\int_{0}^{\infty} s^{p-1}e^{-\underline{k}^1_i s} \int_{0}^{s} |k_i^0 (u+t-s)|^p\mathrm{d} u \mathrm{d} s\bigg\}\nonumber\\
 	:=&\alpha_2(t,p),\\
G_4\le &{11}^{p-1}\max_{1\le i\le n}\bigg\{\bigg(\frac{1}{k_i^-}\bigg)^{\frac{p}{q}}\bigg(\sum_{j=1}^{n}(2CL_j^g)^q\bigg)^{\frac{p}{q}} \int_{-\infty}^t e^{-k_i^-(t-s)} \sum_{j=1}^{n} |d_{ij}^0(s)|^p \mathrm{d} s\bigg\}:=\alpha_3(t,p) ,\\
G_5\le &{11}^{p-1} \max_{1\le i\le n}\bigg\{\bigg(\frac{1}{k_i^-}\bigg)^{\frac{p}{q}}\bigg(\sum_{j=1}^{n} (\bar d_{ij}^1)^q\bigg)^{\frac{p}{q}} \sum_{j=1}^{n}(L_j^g)^p\frac{e^{k_i^-\eta_{ij}^+}}{1-\dot \eta_{ij}^+} \int_{-\infty}^t e^{-k_i^-(t-s)} E|z(s)|_0^p \mathrm{d} s\bigg\},\\
   G_6\leq&{11}^{p-1} \max_{1\le i\le n}\bigg\{\bigg(\frac{1}{\underline{k}_i^1}\bigg)^{\frac{p}{q}}\bigg(\sum_{j=1}^{n} (\bar d_{ij}^1)^q\bigg)^{\frac{p}{q}} \sum_{j=1}^{n}(2CL_j^g)^p\int_{0}^{\infty} s^{p-1}e^{-\underline{k}^1_i s} \int_{0}^{s} |k_i^0 (u+t-s)|^p\mathrm{d} u \mathrm{d} s\bigg\}\nonumber\\
   :=&\alpha_4(t,p),\\
G_7 \le &{11}^{p-1}\max_{1\le i\le n}\bigg\{\bigg(\frac{1}{k_i^-}\bigg)^{\frac{p}{q}} \int_{-\infty}^t e^{-k_i^-(t-s)} |I_i^0(s)|^p \mathrm{d} s\bigg\}:=\alpha_5(t,p) ,\\
 G_8\le& {11}^{p-1} \max_{1\le i\le n}\bigg\{\bigg(\frac{1}{\underline{k}_i^1}\bigg)^{\frac{p}{q}}(\bar{I}_i^1)^p  \int_{0}^{\infty} s^{p-1}e^{-\underline{k}^1_i s}  \int_{0}^{s} |k_i^0 (u+t-s)|^p\mathrm{d} u  \mathrm{d} s \bigg\}
   :=\alpha_6(t,p).\label{433}
\end{align}

When $p>2$, by the BDG inequality and the Fubini theorem, we obtain that
\begin{align}
G_9 \le &{11}^{p-1}\max_{1\le i\le n}\bigg\{C_p\bigg[\int_{-\infty}^{t}\Big(e^{-\int_{s}^{t}k_i(u)\mathrm{d} u}\Big)^{\frac{p}{p-2}} \mathrm{d} s\bigg]^{\frac{p-2}{p}\times \frac{p}{2}}\bigg[\int_{-\infty}^{t}\Big(e^{-\int_{s}^{t}k_i(u)\mathrm{d} u}\Big)^{\frac{p}{2}}\nonumber\\
	&\times E\bigg|\sum_{j=1}^{n}U_{ij}^0(s)\delta_j(x_j(s-r_{ij}(s)))\bigg|^p\mathrm{d} s \bigg] \bigg\}\nonumber\\
	\le &{11}^{p-1}\max_{1\le i\le n}\bigg\{C_p\bigg(\frac{p-2}{pk_i^-}\bigg)^{\frac{p-2}{2}}\bigg(\sum_{j=1}^{n}(2CL_j^\delta)^q\bigg)^{\frac{p}{q}}\int_{-\infty}^{t}e^{-\frac{p}{2}k_i^-(t-s)} \bigg|\sum_{j=1}^{n}U_{ij}^0(s)\bigg|^p\mathrm{d} s \bigg\}\nonumber\\
	:=&\alpha_7(t,p),\\
 G_{10}\le &{11}^{p-1}\max_{1\le i\le n}\bigg\{C_p\bigg[\int_{-\infty}^{t}\Big(e^{-\int_{s}^{t}k_i(u)\mathrm{d} u}\Big)^{ \frac{p}{p-2}} \mathrm{d} s\bigg]^{\frac{p-2}{p}\times \frac{p}{2}}\bigg[\int_{-\infty}^{t}\Big(e^{-\int_{s}^{t}k_i(u)\mathrm{d} u}\Big)^{ \frac{p}{2}}\nonumber\\
    &\times E\bigg|\sum_{j=1}^{n}U_{ij}^1(s)[\delta_j(x_j(s-r_{ij}(s)))-\delta_j(y_j(s-r_{ij}(s)))]\bigg|^p \mathrm{d} s \bigg] \bigg\}\nonumber\\
    \le &{11}^{p-1}\max_{1\le i\le n}\bigg\{C_p\bigg(\frac{p-2}{pk_i^-}\bigg)^{\frac{p-2}{2}}\bigg(\sum_{j=1}^{n} (\bar U_{ij}^1)^q \bigg)^{\frac{p}{q}} \sum_{j=1}^{n}(L_j^\delta)^p\frac{e^{\frac{p}{2}k_i^-r_{ij}^+}}{1-\dot r_{ij}^+} \int_{-\infty}^t e^{-\frac{p}{2}k_i^-(t-s)} E|z(s)|_0^p \mathrm{d} s\bigg\},\\
    G_{11}
\leq&{11}^{p-1}C_p\max_{1\le i\le n}\bigg\{E\bigg(\int_{-\infty}^t e^{-2\underline{k}_i^1 (t-s)} \bigg(\int_s^t |k_i^0 (u)|\mathrm{d} u\bigg)^2  \bigg|\sum_{j=1}^{n}U_{ij}^1(s)\delta_j(y_j(s-r_{ij}(s)))\bigg|^2\mathrm{d} s\bigg)^{\frac{p}{2}}\bigg\}\nonumber\\
\leq&{11}^{p-1}C_p\max_{1\le i\le n}\bigg\{\bigg(\frac{p-2}{p\underline{k}_i}\bigg)^{\frac{p-2}{2}}\bigg(\sum_{j=1}^{n} (\bar U_{ij}^1)^q\bigg)^{\frac{p}{q}} \sum_{j=1}^{n}(2CL_j^\delta)^p\nonumber\\
	 &\times \int_{-\infty}^t (t-s)^{p-1}e^{-\frac{p}{2}\underline{k}_i^1 (t-s)} \int_s^t |k_i^0 (u)|^p\mathrm{d} u  \mathrm{d} s \bigg\}\nonumber\\
\leq&{11}^{p-1}C_p\max_{1\le i\le n}\bigg\{\bigg(\frac{p-2}{p\underline{k}_i}\bigg)^{\frac{p-2}{2}}\bigg(\sum_{j=1}^{n} (\bar U_{ij}^1)^q\bigg)^{\frac{p}{q}} \sum_{j=1}^{n}(2CL_j^\delta)^p\nonumber\\
	 &\times \int_{0}^\infty s^{p-1}e^{-\frac{p}{2}\underline{k}_i^1 s} \int_{0}^s |k_i^0 (u+t-s)|^p\mathrm{d} u  \mathrm{d} s \bigg\}
    :=\alpha_8(t,p).\label{436}
\end{align}
	When $p=2$, according to the It\^{o} isometry and the Fubini theorem
\begin{align}
G'_9\le &11\max_{1\le i\le n}\bigg\{\sum_{j=1}^{n}(2CL_j^\delta)^2 \int_{-\infty}^t e^{-2k_i^- (t-s)}\bigg|\sum_{j=1}^{n}U_{ij}^0(s)  \bigg|^2\mathrm{d} s \bigg\}:=\tilde \alpha_7(t),\label{437}\\
G'_{10}\le &11\max_{1\le i\le n}\bigg\{\sum_{j=1}^{n}(\bar U_{ij}^1)^2\sum_{j=1}^{n}(L_j^\delta)^2\frac{e^{2k_i^- r_{ij}^+}}{1-\dot r_{ij}^+}\bigg\}\int_{-\infty}^t e^{-2k^- (t-s)} E|
z_0(s)|^2 \mathrm{d} s,\\
G'_{11}\le& {11}^{p-1} \max_{1\le i\le n}\bigg\{\sum_{j=1}^{n} (\bar U_{ij}^1)^2 \sum_{j=1}^{n}(2CL_j^\delta)^2 \int_{0}^{\infty} s e^{-2\underline{k}^1_i s}\int_{0}^{s}|k_i^0 (u+t-s)|^2\mathrm{d} u \mathrm{d} s \bigg\}
	:=\tilde \alpha_8(t).\label{4.67}
  \end{align}

For $p>2$,  from \eqref{3.58}-\eqref{436}, it follows that
   \begin{align*}
   	E|z(t)|_0^p \le \alpha(t,p)+\beta \int_{-\infty}^{t} e^{-k^-(t-s)} E|z(s)|_0^p \mathrm{d} s,
   \end{align*}
where $\alpha(t):=\sum_{\varsigma=1}^{8}\alpha_\varsigma(t,p)$, $k^-=\min_{1\le i\le n}\{k_i^-\}$ and $\beta$ is defined in $(H_6)$.

In view of Lemma \ref{G}, for any $\gamma \in (0,k^- -\beta)$, there holds
 \begin{align}\label{4.68}
	E|z(t)|_0^p \le \alpha(t,p)+\beta \int_{-\infty}^{t} e^{-\gamma (t-s)} \alpha(s,p) \mathrm{d} s.
\end{align}
Note that
\begin{align*}
&\lim_{h\rightarrow \infty} \frac{1}{\nu ([-h,h])} \int_{-h}^{h}\int_{-\infty}^{t}e^{-k_i^-(t-s)}|c_{ij}^0(s)|^p \mathrm{d} s \mathrm{d} \mu (t)\\
\leq&\int_{0}^{+\infty}e^{-k_i^-  s }\lim_{h\rightarrow \infty} \frac{1}{\nu ([-h,h])} \int_{-h}^{h} |c_{ij}^0(t-s)|^p  \mathrm{d} \mu (t) \mathrm{d} s\\
\le &\int_{0}^{+\infty}e^{-k_i^-  s }\lim_{h\rightarrow \infty} \frac{\nu_s ([-h-s,h+s])}{\nu ([-h,h])} \times \frac{1}{\nu_s ([-h-s,h+s])} \int_{-h-s}^{h+s} |c_{ij}^0(t)|^p  \mathrm{d} \mu_s (t) \mathrm{d} s=0,
\end{align*}
\begin{align*}
&\lim_{h\rightarrow \infty} \frac{1}{\nu ([-h,h])} \int_{-h}^{h}\int_{0}^{\infty} s^{p-1}e^{-\underline{k}^1_i s} \int_{0}^{s} |k_i^0 (u+t-s)|^p\mathrm{d} u \mathrm{d} s \mathrm{d} \mu (t)\\
\leq&\int_{0}^{\infty} s^{p-1}e^{-\underline{k}^1_i s} \int_{0}^{s}\lim_{h\rightarrow \infty} \frac{1}{\nu ([-h,h])} \int_{-h}^{h} |k_i^0 (u+t-s)|^p\mathrm{d} \mu (t)\mathrm{d} u \mathrm{d} s=0,
\end{align*}

\begin{align*}
&\lim_{h\rightarrow \infty} \frac{1}{\nu ([-h,h])} \int_{-h}^{h}\int_{-\infty}^{t}e^{-\gamma (t-s)}\int_{-\infty}^{s} e^{-k_i^-(s-u)} |c_{ij}^0(u)|^p \mathrm{d} u \mathrm{d} s \mathrm{d} \mu (t) \nonumber\\
\leq&\int_{0}^{+\infty}e^{-\gamma s}\lim_{h\rightarrow \infty} \frac{1}{\nu ([-h,h])} \int_{-h}^{h}\int_{-\infty}^{t-s} e^{-k_i^-(t-s-u)} |c_{ij}^0(u)|^p \mathrm{d} u \mathrm{d} \mu (t) \mathrm{d} s \nonumber\\
\leq&\int_{0}^{+\infty}e^{-\gamma s} \int_{0}^{\infty} e^{-k_i^- u}\lim_{h\rightarrow \infty}\frac{1}{\nu ([-h,h])} \int_{-h}^{h} |c_{ij}^0(t-s -u)|^p \mathrm{d} \mu (t)\mathrm{d} u \mathrm{d} s=0
\end{align*}
and
\begin{align*}
&\lim_{h\rightarrow \infty} \frac{1}{\nu ([-h,h])} \int_{-h}^{h}\int_{-\infty}^{t}e^{-\gamma (t-s)}  \int_{0}^{\infty} v^{p-1}e^{-\underline{k}^1_i v} \int_{0}^{v} |k_i^0 (u+s-v)|^p\mathrm{d} u \mathrm{d} v \mathrm{d} s\mathrm{d}\mu(t)\\
\leq&\int_{0}^{\infty}e^{-\gamma s} \int_{0}^{\infty} v^{p-1}e^{-\underline{k}^1_i v} \int_{0}^{v}\lim_{h\rightarrow \infty} \frac{1}{\nu ([-h,h])} \int_{-h}^{h} |k_i^0 (u+t-s-v)|^p \mathrm{d}\mu(t)\mathrm{d} u\mathrm{d} v \mathrm{d} s=0.
\end{align*}

We can get, for $j=1,2$,
\begin{align*}
	\lim_{h\rightarrow \infty} \frac{1}{\nu ([-h,h])} \int_{-h}^{h} \alpha_j(t,p) \mathrm{d} \mu (t)=	\lim_{h\rightarrow \infty} \frac{1}{\nu ([-h,h])} \int_{-h}^{h}\int_{-\infty}^{t}e^{-\gamma (t-s)} \alpha_j(s,p) \mathrm{d} s \mathrm{d} \mu (t)=0.
\end{align*}

Analogously, for $i=3,4,\ldots,8$, one can get
\begin{align*}
	\lim_{h\rightarrow \infty} \frac{1}{\nu ([-h,h])} \int_{-h}^{h} \alpha_i(t,p) \mathrm{d} \mu (t)=	\lim_{h\rightarrow \infty} \frac{1}{\nu ([-h,h])} \int_{-h}^{h}\int_{-\infty}^{t}e^{-\gamma (t-s)} \alpha_i(s,p) \mathrm{d} s \mathrm{d} \mu (t)=0.
\end{align*}

Consequently, by \eqref{4.68},  we arrive at
\begin{align}\label{fn1}
\lim_{h\rightarrow \infty} \frac{1}{\nu ([-h,h])} \int_{-h}^{h}E|z(t)|_0^p \mathrm{d} \mu (t) =0.
\end{align}

For $p=2$,  by \eqref{3.58}-\eqref{433}, \eqref{437}-\eqref{4.67} and Lemma \ref{G}, for any $\tilde \gamma \in (0,k^- -\tilde \beta)$, it holds
 \begin{align}\label{4.69}
 	E|z(t)|_0^2 \le \tilde \alpha(t)+ \tilde \beta \int_{-\infty}^{t} e^{- \tilde \gamma (t-s)} \tilde \alpha(s) \mathrm{d} s,
\end{align}
where $\tilde \alpha(t):=\sum_{\nu=1,2,\ldots,6}\alpha_\nu(t,2)+\tilde \alpha_7(t)+\tilde \alpha_8(t)$ and $\tilde{\beta}$ is given in $(H_6)$.

Using \eqref{4.69}, similar to the proof of the case for $p>2$, we obtain that
\begin{align}\label{fn2}
	\lim_{h\rightarrow \infty} \frac{1}{\nu ([-h,h])} \int_{-h}^{h}E|z(t)|_0^2 \mathrm{d} \mu (t) =0.
\end{align}

According to Theorem \ref{key}, \eqref{fn1} and \eqref{fn2}, we conclude that $z\in \mathcal{O}_D^p(\mathbb{R},
	\mathcal{L}^p(\Omega,\mathbb{X}),\mu,\nu)$. The proof is finished.
\end{proof}

By a straightforward modification of the argument employed in the proof of Theorem 3.7 in \cite{liw}, the following theorem can be obtained.
\begin{theorem}\label{thm5}
	Assume that $(H_1)-(H_4)$ and $(H_6)$ hold. If $x$ is Besicovitch  $(\mu,\nu )$-pseudo almost periodic solution in distribution with initial value $\phi$ and $x^*$ is an arbitrary solution of system \eqref{e1} with initial value $\phi^*$, then there exist positive numbers $\lambda>0$ and $M>0$ satisfying
	\begin{align*}
		E|x(t)-x^*(t)|_0^p \le ME|\phi - \phi^*|_\rho^pe^{-p\lambda t}, \,\, t>0,
	\end{align*}
	where
	$E|\phi - \phi^*|_\rho^p =\max_{1\le i\le n}\big\{\sup_{s \in [-\rho,0] }E|\phi_i (s)-a_i(t)\phi^*_i (s) |^p\big\}.$
\end{theorem}

\section{An example}
\setcounter{equation}{0}
\indent

In system \eqref{e1}, let $n=2$, and take
\begin{align*}
	f_j(x_j) =& 0.05 \sin (x_j),\,
	g_j(x_j) =0.04\tanh (x_j), \,
	\delta_j (x_j) = 0.05\sin (x_j), j=1,2,\\
	k^1_1(t)=&1.2-0.2\cos \sqrt 2t -0.1e^{-|t|},\,
	k^1_2(t)=1.1-0.1\sin \sqrt 2t -0.05\frac{\sin t}{1+t^2},\\
	c_{11}^1(t)=&0.02+0.008\bigg(\cos t+\frac{\sin t}{1+t^2}\bigg),\,
	c_{12}^1(t)= 0.03 + 0.002\bigg(\sin \sqrt 5t +\frac{1}{1+t^2}\bigg),\\
	c_{21}^1(t)=&0.04 + 0.001(\cos \sqrt 2t+e^{-|t|}),\,
	c_{22}^1(t)=0.03 + 0.004\sin \sqrt 3t,\\
	d_{11}^1(t)=&0.02+0.007(\cos t+e^{-|t|}),\,
	d_{12}^1(t)=0.02+0.001(\sin \sqrt 2t+e^{-|t|}),\\
	d_{21}^1(t)=&0.03+0.002\bigg(\sin \sqrt 2t+\frac{1}{1+t^2}\bigg),
	d_{22}^1(t)=0.04+0.001\sin \sqrt 2t,\\
	U_{11}^1(t)=&0.03+0.001(\cos \sqrt 2t+e^{-|t|}),\,
	U_{12}^1(t)=0.02+0.008\sin \sqrt 2t,\\
	U_{21}^1(t)=&0.03+0.002 (\sin \sqrt 2t+e^{-|t|} ),\,
	U_{22}^1(t)=0.04+0.001\bigg(\cos \sqrt 2t+\frac{\sin t}{1+t^2}\bigg),\\
	I^1_1(t)=&0.8\cos\sqrt{2} t+\frac{\sin t}{1+t^2}+0.2\sin \sqrt{3} t,\,
   I^1_2(t)=0.6 \sin \sqrt{3} t+0.2\cos t+e^{-|t|},\\
   k_1(t)=&k_1^1(t)-0.1e^{-|t|},\,
	k_2(t)=k_2^1(t)-\frac{0.05\sin t}{1+t^2}, I_1(t)=I_1^1(t),\,
   I_2(t)=I_2^1(t)+\mathrm{sign}(t),\\
    c_{11}(t)=&c_{11}^1(t)+O_1(t), c_{12}(t)=c_{12}^1(t)+O_2(t),\,
	c_{21}(t)=c_{21}^1(t)+O_2(t), c_{22}(t)=c_{22}^1(t),\\
	d_{11}(t)=&d_{11}^1(t)+O_1(t), d_{12}(t)=d_{12}^1(t)+O_2(t),\,
	d_{21}(t)=d_{21}^1(t)-O_2(t), d_{22}(t)=d_{22}^1(t),\\
	U_{11}(t)=&U_{11}^1(t)+O_2(t), U_{12}(t)=U_{12}^1(t),\,
	U_{21}(t)=U_{21}^1(t)+O_1(t), U_{22}(t)=U_{22}^1(t),\\
 \sigma_1(t)=&0.04+0.02\sin t, \sigma_2(t)=0.05+0.01\cos t,\,
   \eta_{11}(t)=0.03+0.02\cos \sqrt{2} t,\\
    \eta_{12}(t)=&0.04+0.02\cos t,\,
    \eta_{21}(t)=0.05+0.01\sin t,\,
	\eta_{22}(t)=0.05+0.01\sin \sqrt{2} t,\\
	r_{11}(t)=&0.03+0.03\cos\sqrt{3} t,\,
	r_{12}(t)=0.04+0.02\sin t,\\
	r_{21}(t)=&0.03+0.02\cos t,\,
	r_{22}(t)=0.04+0.02\sin t,
\end{align*}
where $\mathrm{sign}(t)$ indicate the sign function, and
\begin{align*}
	O_1(t) = \begin{cases}
		0.004e^{t},&t \le 0,\\
		0.004,&t > 0
	\end{cases}\quad \mathrm{and}  \quad
	O_2(t) = \begin{cases}
		0.002e^{t},&t \le 0,\\
		0.002\cos t,&t > 0.
	\end{cases}
\end{align*}
Obviously, $|O_1|_{B^p}\ne 0$ and $|O_2|_{B^p}\ne 0$.

Take the Radon-Nikodym derivatives of $\mu$ and $\nu$ are $\sin 2 t+2$ and $t^2+\cos t+4$ with respect to the Lebesgue measure on $\mathbb{R}$, respectively. Then, we have
\begin{align}\label{zv}
 \overline{\lim_{h\to \infty}} \frac{\mu([-h,h])}{\nu([-h,h])}=\overline{\lim_{h\to \infty}} \frac{\int_{-h}^{h}(\sin 2t+2)\mathrm{d}t}{\int_{-h}^{h}(t^2+\cos t+4)\mathrm{d}t}=0<\infty,
\end{align}
and hence, condition $(A_1)$ is satisfied. Moreover, for all $\sigma\in \mathbb{R}$ and $A\in \Sigma$, we have
$\mu_\sigma(A)\leq 3 \mu (A).$
Take $\zeta=3$ and $\mathcal{J}=\emptyset$, then  condition  $(A_2)$ is also satisfied.

Since the Besicovitch seminorns of function sign$(t),e^{-|t|},\frac{\sin t}{1+t^2}, \frac{1}{1+t^2}$ are $0$, these functions are Besicovitch almost periodic.
Since $O_1(t), O_2(t), \mathrm{sign}(t)$ are bounded functions, in view of \eqref{zv}, $O_1, O_2, \mathrm{sign}(\cdot) \in\mathcal{O}^p (\mathbb{R},\mathbb{R},\mu,\nu)$. Therefore, condition $(H_1)$ is satisfied.
Also, it is clear that $(H_2)$ is fulfilled.

Through computation, we have
\begin{gather*}
  L_j^f =0.05, L_j^g=0.04, L_j^\delta =0.05, k_1^+=1.4, k_2^+=1.2, k_1^-=0.800, k_2^-=0.9590, k^-=0.800, \\
\underline{k}_1^1=0.9,\, \underline{k}_2^1=0.9797, \,  \bar c_{11}^1=0.0310,\, \bar c_{12}^1=0.03342,\,  \bar c_{21}^1=0.04200,\, \bar c_{22}^1=0.03400,\\
	 \bar d_{11}^1=0.03400,\bar d_{12}^1=0.02136, \bar d_{21}^1=0.03301,\bar d_{22}^1=0.04100, \bar U_{11}^1=0.03200, \bar U_{12}^1=0.02800,\\
	 \bar U_{21}^1=0.03273,\bar U_{22}^1=0.04119, c_{11}^+=0.03410, c_{12}^+=0.03516, c_{21}^+=0.04400, c_{22}^+=0.03400,\\
	 d_{11}^+=0.03800, d_{12}^+=0.02306, d_{21}^+=0.03218,d_{22}^+=0.04100, U_{11}^+=0.03400, U_{12}^+=0.02800,\\
	U_{21}^+=0.03673,\, U_{22}^+=0.04119,\,\, \sigma_1^+=0.06, \,\,\sigma_2^+=0.06,\,\, \eta_{11}^+=0.05, \,\,\eta_{12}^+=0.06, \,\,\eta_{21}^+=0.06, \\
	\eta_{22}^+=0.06, r_{11}^+=0.06, r_{12}^+=0.06, r_{21}^+=0.05, r_{22}^+=0.06,\dot \sigma_1^+=0.02, \dot \sigma_2^+=0.01,\\
 \dot \eta_{11}^+=\dot \eta_{12}^+=0.02,\quad \dot \eta_{21}^+=\dot \eta_{22}^+= 0.01,\quad \dot r_{11}^+= 0.03, \quad\dot r_{12}^+=\dot r_{21}^+=\dot r_{22}^+=0.02,
\end{gather*}
and when $p=5$, 	$F_1 \approx  3.058\times 10^{-7}$, when $p=2$, $ F_2 \approx  1.0505\times 10^{-4}$,
when $p=5$, 	$F'_1 \approx  2.24\times 10^{-7}$, when $p=2$, $ F'_2 \approx  1.10\times 10^{-4}$,
when $p=5$,   	$\Xi_1 \approx 1.9164\times 10^{-4}<1 $,  when $p=2$, $\Xi_2 \approx 5.61636\times10^{-4}<1$,
when $p=5$,   	$\beta \approx 3.03\times10^{-5}<k^-=0.8$, when $p=2$, $\tilde{\beta} \approx 4.489\times10^{-4}<k^-=0.8$.
Hence, $(H_1)$, $(H_4)$ and $(H_6)$ are satisfied.

By Theorem \ref{thm5}, the system \eqref{e1} has a unique Besicovitch $(\mu,\nu)$-pseudo almost periodic solution in distribution that is global exponential stability (see Figures 1-2).
	\begin{figure}[!h]
	\centering
	\includegraphics[width=16cm,height=7cm]{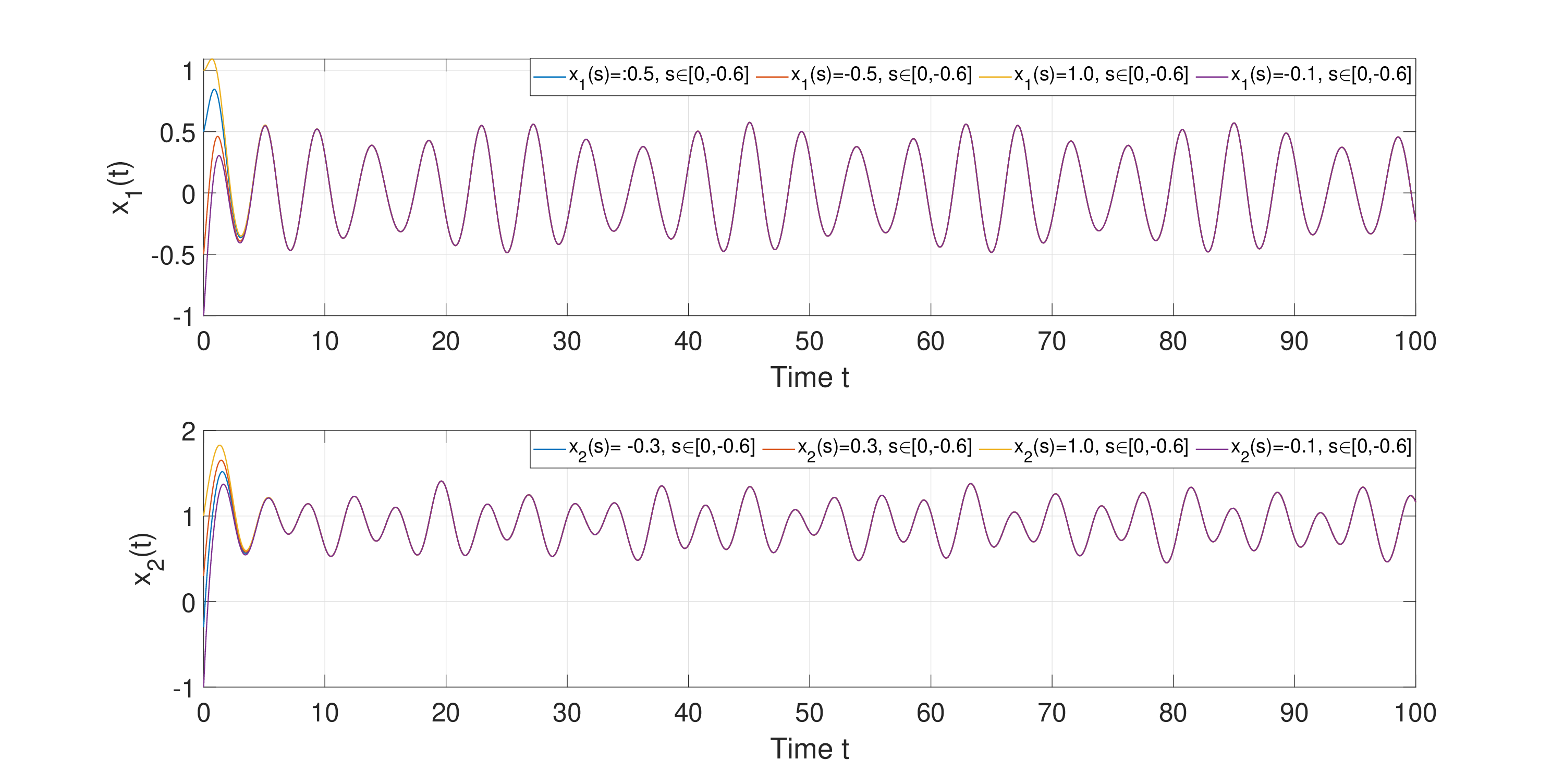}
	\caption{States $x_1$ and $x_2$ of \eqref{e1} with different initial values.}
\end{figure}
	\begin{figure}[!h]
	\centering
	\includegraphics[width=16cm,height=6cm]{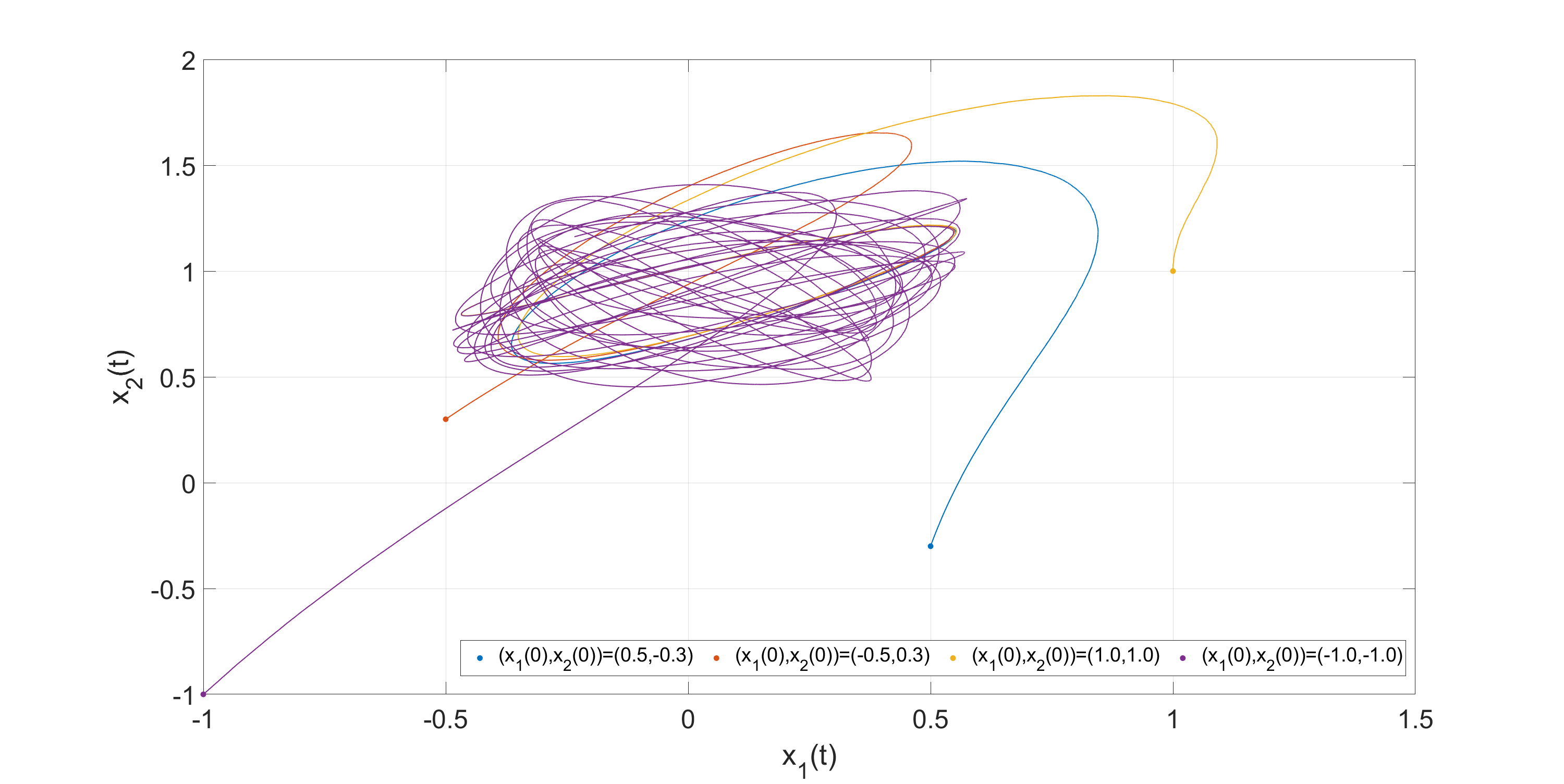}
	\caption{Phase portraits (circles mark initial points).}
\end{figure}

\begin{remark}
In the example, about the decomposition of $I_2(t)$, we have taken $I_2(t)=I^1_2(t)+\mathrm{sign}(t)$, where  $I^1_2(t)=0.6 \sin \sqrt{3} t+0.2\cos t+e^{-|t|}$. In fact, we can also decompose it as $I_2(t)=\breve{I}^1_2(t)+e^{-|t|}+\mathrm{sign}(t)$, where  $\breve{I}^1_2(t)=0.6 \sin \sqrt{3} t+0.2\cos t$. This shows that the decomposition in Definition \ref{def333} may  not unique.
\end{remark}

\section{Conclusion}
\indent	

In this paper, we have proposed the concept of pseudo almost periodic stochastic processes in the Besicovitch sense, thereby filling the gap that the concept of Besicovitch pseudo almost periodic stochastic processes has not yet been systematically established. As an application of this concept, we have studied the existence of Besicovitch pseudo almost periodic solutions in distribution for a class of stochastic delayed NNs in which all parameters except the time delays are pseudo almost periodic and established stability criteria for such solutions in the $p$-th mean sense. It should be emphasized that stochastic NNs serve only as an application vehicle for the new concept proposed in this paper; the core contribution of this paper lies in providing a more general descriptive tool for the complex oscillatory behavior of stochastic non-autonomous systems than the existing  the Besicovitch almost periodic framework.

The results of this paper show that combining the pseudo almost periodic structure with the Besicovitch space, handling the existence of solutions in the distribution sense and stability in the $p$-th mean sense, is a feasible and natural analytical path. This path not only retains the advantages of the Besicovitch almost periodic concept, namely broad coverage and relatively weak conditions, but also accommodates a wider class of stochastic non-autonomous perturbations through the pseudo almost periodic decomposition. However, the work of this paper remains preliminary:  the application part only examines one class of stochastic NNs with delays. Future work can proceed in the following directions: first, studying Besicovitch pseudo almost periodic solutions under more general stochastic terms (such as Markov switching and L\'{e}vy noise); second, establishing criteria for stability in distribution, almost sure stability, and other stability notions different from the stability in the $p$-th mean sense; third, providing numerical examples to verify the operability of the theoretical results; and fourth, extending the concept to broader systems such as stochastic partial differential equations and abstract stochastic functional differential equations.

\section*{Data availability}
\indent

 No datasets were generated or analysed during the current study.

\section*{Conflicts of Interest}
\indent

The author declares that there are no conflicts of interest.

{}


\begin{thebibliography}{}

\bibitem{bo} H. Bohr, Zur Theorie der fastperiodischen Funktionen I, Acta Math. 45 (1925) 29–127.
	
\bibitem{st} V. V. Stepanov, Sur quelques g\'en\'eralisations des fonctions presque p\'eriodiques, CR Acad. Sci. Paris 181 (1925) 90–92.
	
\bibitem{we} H. Weyl, Integralgleichungen und fastperiodische Funktionen, Math. Ann. 97 (1927) 338–356.
	
\bibitem{be}A. S. Besicovitch, On generalized almost periodic functions, Proc. Lond. Math. Soc., 25 (1926), 495-512.

\bibitem{y1} J. Andres, A.M. Bersani, R.F. Grande, Hierarchy of almost-periodic function spaces, Rend. di Mat. 26 (2006) 121–188.

\bibitem{wave} C. Corduneanu, Almost periodic oscillations and waves, Springer, New York, (2009).

\bibitem{wave1}T. Diagana, Almost automorphic type and almost periodic type functions in abstract spaces,  Springer,  New York, (2013).

\bibitem{y2} Y. Li, X. Huang, Besicovitch almost periodic solutions to stochastic dynamic equations with delays, Qual. Theory Dyn. Syst. 21 (2022) 74.

%Y. Li, X. Wang, Besicovitch almost periodic stochastic processes and almost periodic solutions of Clifford-valued stochastic neural networks, Discrete Contin. Dyn. Syst.  28(3) (2023) 2154-2183.

\bibitem{y3}P. H. Bezandry, Existence of almost periodic solutions for semilinear stochastic evolution equations driven by fractional Brownian motion, Electron, J. Differ. Equ. 2012 (2012), no. 156, 1–21.

\bibitem{y5}Y. Sun, J. Cao, $p$th moment exponential stability of stochastic recurrent neural networks with time-varying delays, Nonlinear Anal.: Real World Appl. 8(4) (2007) 1171-1185.

\bibitem{y6} J. Pan, G. Chen, Generalized Halanay inequalities and $p$th moment exponential stability of impulsive neutral stochastic neural networks with time varying delays, Quaest. Math. 2026, in press, https://doi.org/10.2989/16073606.2026.2750939

\bibitem{y4} Y. Li, M. Huang, B. Li,  Besicovitch almost periodic solutions for fractional-order quaternion-valued neural networks with discrete and distributed delays, Math. Meth.  Appl. Sci. 45 (2022) 4791-4808.
	


\bibitem{2026}Y. Li, Measure pseudo Besicovitch almost periodic functions and an application to a Lasota-Wazewska model, Qual. Theorey Dyn. Syst., 25 (2026), 122.

\bibitem{W} Y. Li, X. Wang, Besicovitch almost periodic stochastic processes and almost periodic solutions of Clifford-valued stochastic neural networks, Discret. Contin. Dyn. Syst. B, 28 (2023), 2154–2183.

\bibitem{measure1}J. Blot, P. Cieutat, K. Ezzinbi, Measure theory and pseudo almost automorphic functions: New developments and applications, Nonlinear Anal. 75 (2012) 2426-2447.

\bibitem{a3} S. Hu, C. Huang, F. Wu, Stochastic differential equations, Science Press, Beijing, (2008).

\bibitem{27}M. Kamenskii, O. Mellah, P. Fitte, Weak averaging of semilinear stochastic differential equations with almost periodic coefficients, J. Math. Anal. Appl. 427 (2015)  336-364.



\bibitem{liw} Y. Li, X. Li, Besicovitch almost periodic stochastic processes and almost periodic solutions of Clifford-valued stochastic neural networks, Discrete Contin. Dyn. Syst.-B 28(3) (2023) 2154-2183.
	
	
	
	
	
		
		


	\end{thebibliography}
\end{document}